\documentclass[12pt]{amsart}
\usepackage{amsmath,amssymb}
\usepackage{graphicx,epsfig,subfigure,psfrag}
\begin{document}

%
%
\newtheorem{theorem}{Theorem}
\newtheorem{proposition}[theorem]{Proposition}
\newtheorem{lemma}[theorem]{Lemma}
\newtheorem{corollary}[theorem]{Corollary}
\newtheorem{definition}[theorem]{Definition}
\newtheorem{remark}[theorem]{Remark}
\numberwithin{equation}{section}
\numberwithin{theorem}{section}
\newcommand{\be}{\begin{equation}}
\newcommand{\ee}{\end{equation}}
\newcommand{\re}{\mathbb{R}}
\newcommand{\n}{\nabla}
\newcommand{\ren}{\mathbb{R}^N}
\newcommand{\iy}{\infty}
\newcommand{\pa}{\partial}
\newcommand{\ms}{\medskip\vskip-.1cm}
\newcommand{\mpb}{\medskip}
\newcommand{\BB}{{\bf B}}
\newcommand{\Am}{{\bf A}_{2m}}
\newcommand{\bL}{\BB^*}
\newcommand{\bLs}{\BB}
\renewcommand{\a}{\alpha}
\renewcommand{\b}{\beta}
\newcommand{\g}{\gamma}
\newcommand{\ka}{\kappa}
\newcommand{\G}{\Gamma}
\renewcommand{\d}{\delta}
\newcommand{\D}{\Delta}
\newcommand{\e}{\varepsilon}
\renewcommand{\l}{\lambda}
\renewcommand{\o}{\omega}
\renewcommand{\O}{\Omega}
\newcommand{\s}{\sigma}
\renewcommand{\t}{\tau}
\renewcommand{\th}{\theta}
\newcommand{\z}{\zeta}
\newcommand{\wx}{\widetilde x}
\newcommand{\wt}{\widetilde t}
\newcommand{\noi}{\noindent}
\newcommand{\lb}{\left (}
\newcommand{\rb}{\right )}
\newcommand{\lsb}{\left [}
\newcommand{\rsb}{\right ]}
\newcommand{\lab}{\left \langle}
\newcommand{\rab}{\right \rangle }
\newcommand{\gap}{\vskip .5cm}
\newcommand{\bz}{\bar{z}}
\newcommand{\bg}{\bar{g}}
\newcommand{\Ba}{\bar{a}}
\newcommand{\bt}{\bar{\th}}
\def\com#1{\fbox{\parbox{6in}{\texttt{#1}}}}

\title{\bf On centre subspace behaviour \\ in thin film equations}

\author
{V A Galaktionov $\dag$ and P J Harwin$\ddag$}

\address{$\dag$ Department of Math. Sci., University of Bath,
 Bath, BA2 7AY}
\email{vag@maths.bath.ac.uk}

\address{$\ddag$ Department of Math. Sci., University of Bath,
  Bath, BA2 7AY, UK}
  \email{P.J.Harwin@maths.bath.ac.uk}

\keywords{Quasilinear thin  film equation, critical absorption
exponent, similarity solutions, asymptotic behaviour. To appear in SIAM J. Appl. Math.}
\date{\today}



\begin{abstract}

The large-time behaviour of weak nonnegative solutions of the thin
film equation (TFE) with absorption
\[
u_t = -\nabla \cdot (|u|^n \nabla \Delta u) - |u|^{p-1}u, 
\]
with parameters $n \in (0,3)$  and
  $p>1$, is studied.
The standard free-boundary problem with zero-height, zero contact
angle, and zero-flux conditions at the interface
   and bounded
compactly supported initial data is considered.  It is shown that
there exists the {\em critical} absorption exponent
\[
\textstyle{p_0= 1+n + \frac 4 N}
\]
such that, for $p=p_0$, the asymptotic behaviour of solutions
$u(x,t)$ for $t \gg 1$ is represented by the  well-known
source-type solution of the pure TFE  absorption,
\[
\textstyle{u_s(x,t) = t^{-  \b N} F(y), \quad y= x/t^ \b, \quad
\mbox{with the exponent} \,\,\, \b = \frac 1{4+nN},}
\]
which is perturbed by a couple of $\ln t$-factors.
For $n=1$, this behaviour is associated with  the centre subspace
for the rescaled linearized thin film operator and is given by
\[
\textstyle{u(x,t) \sim (t \ln t)^{- \b N} F( x/t^\b (\ln t)^{- \b
N /4}), \quad \mbox{with} \,\,\, \b = \frac 1{4+N},}
\]
where $F(y)= \frac 1{8(N+2)(N+4)}{(a_*^2-|y|^2)^2}$ and the
constant  $a_*>0$ depends on dimension $N$ only.
The $2m$th-order generalization of such TFEs with critical
absorption is considered  and
 some local and asymptotic features of changing sign similarity
solutions of the Cauchy
problem are described. 

 Our study 
is motivated by  the phenomenon of logarithmically  perturbed
source-type behaviour for the second-order porous medium equation  
with critical absorption
\[
\textstyle{ u_t = \n \cdot (u^n \n u) - u^p  \quad \mbox{in}
\,\,\, \ren \times \re_+, \quad p_0= 1+n + \frac 2 N,} \,\,\, n
\ge 0,
\]
which has been
  known since the 1980s.

\end{abstract}

\maketitle

\section{\sc Introduction: The model, motivation, and results}
\label{Sect1}

Our goal is to describe some  unusual  
asymptotic phenomena for {\em higher-order  quasilinear}
degenerate parabolic equations, in which the nonlinear interaction
between operators involved  deforms  the scaling-invariant
structure of solutions for large times. These  delicate cases of
asymptotic phenomena, such as {\em logarithmic perturbations} of
fundamental or source-type solutions, have been known since the
1980s
 for quasilinear second-order reaction-diffusion equations.
 For {\em semilinear} higher-order parabolic equations, those phenomena
 can be
 detected by using spectral theory of non self-adjoint operators and
semigroup  approaches. For {\em quasilinear} models, similar
asymptotic patterns were unknown.

In the present paper, we introduce a new quasilinear parabolic
model by adding to the standard thin film operator an extra
absorption term. This creates a non-conservative evolution PDE,
  which enjoys a variety of
logarithmically perturbed non-scaling asymptotics in both
free-boundary and the Cauchy problem.
  We then fix several  similarities with
simpler second-order diffusion-absorption models.

We begin with some physical motivation of such models.

\subsection{On general thin film models: a class of conservative and non-conservative PDEs}

For a long time, modern thin film theory and application dealt
with rather complicated nonlinear models. Typically, such models
include the principal quasilinear fourth-order operator and
several
 lower-order 
  terms. For instance,
  the {\em Benney
equation} (1966)
 describes  the nonlinear dynamics of the interface of  2D liquid
films flowing on a fixed inclined plane \cite{Benney},
 \be
 \label{Ben1}
 \mbox{$
 u_t + \frac {2{\rm Re}}3  (u^3)_x + \e \bigl[\bigl( \frac {8 {\rm Re}^2}{15} \, u^6 -
 \frac {2 {\rm Re}}3  \cot \theta \, u^3 \bigr)u_x + \Sigma \, u^3 u_{xxx}\bigr]_x=0,
  $}
  \ee
 where ${\rm Re}$ is the unit-order {\em Reynolds number} of the flow driven by
 gravity, $\s$ is the rescaled {\em Weber number} (related to surface
 tension $\s$), $\theta$ is the angle of plane inclination to the
 horizontal, and $\e= \frac d \l \ll 1$, with $d$ being the average
 thickness of the film and $\l$  the wavelength of the
 characteristic  interfacial disturbances. See \cite{Oron02}.

 TFEs can include
 non-power nonlinearities. For instance, in the  multi-dimensional geometry,
  a typical example is
   \be
   \label{cap11}
    \mbox{$
u_t + \n \cdot \bigl[\bigl(-G \, u^3 + \frac{B M \, u^2}{2P(1+B \,
u)^2} \bigr) \n u \bigr] + S \, \n \cdot (u^3 \n \Delta  u) =0
 $}
 \ee
that describes, in the dimensionless form, the dynamics of a film
in $\re^3$ subject to the actions
 of  thermocapillary, capillary, and gravity forces. Here, $G$,
 $M$, $P$, $B$, and $S$ are the gravity, Marangoni, Prandtl, Biot, and
 inverse capillary numbers respectively. On Marangoni
 instability in such TFE models, see \cite{Oron00}.

The above conservative PDEs preserve the finite mass of thin
films. Non-conservative TFEs occur for evaporating/condensing
films and via other effects, \cite{Oron97, Govor05}.
 Actually, the first study of the vapor thrust  effects
 in the Rayleigh--Taylor instability of an evaporating
 liquid-vapor interface above a hot horizontal wall
was performed by Bankoff in 1961. His stability analysis in 1971
of an evaporating thin liquid film on a hot inclined wall extended
earlier results of Yih (1955, 1963) and Benjamin (1957). The
history and detailed  derivation of  models of (a) evaporating
thin film and (b)
 a condensing thin film,
  can be found in \cite[pp.~946--949]{Oron97}.
A typical TFE of that type in 1D is as follows
\cite[p.~949]{Oron97}:
 \be
 \label{oron1}
 \mbox{$
  u_t + \frac{\bar E}{u+K} + \frac 13 \frac 1{\bar
 C}\big(u^3 u_{xxx}\big)_x
+ \, \big\{\big[ \frac A u +\frac{ \bar E^2}
 D \big(\frac {u}{u+K}\big)^3 + \frac {K M}{\rm Pr}
 \big(\frac {u}{u+K}\big)^2\big]
 u_x\big\}_x =0.
  $}
  \ee
Here, the six terms represent, respectively, the rate of
volumetric accumulation, the mass loss, the stabilization
capillary, van der Waals, vapor thrust, and thermocapillary
 effects.
In the second absorption-like term, $\bar E$ is the scaled
evaporation number and $K$ is the scaled intefacial thermal
resistance that physically represents  a temperature jump from the
liquid surface temperature to the uniform temperature of the
saturated vapor. $D$ is a unit-order scaled ratio between the
vapor and liquid densities.

Another origin of non-conservative TFEs with more complicated
non-divergent operators is the study of flows on a rotating disc
(centrifugal spinning as an efficient mean of coating planar
solids with thin films).
 This gives extra absorption-like, spatially non-autonomous terms
in the equations written in radial geometry, e.g.,
\cite[p.~955]{Oron97}
 \be
 \label{oron2}
  \begin{matrix}
 u_t + \frac 23 E + \frac 1{3r}\big[r^2 u^3 + \e{\rm Re}\big(\frac
 5{12}E r^2 u^4- \frac{34}{105}r^2 u^7\big)\big] \smallskip\smallskip\\
  + \,\frac \e 3 \big\{{\rm Re}\big(\frac 25 r^3u^6-r \frac 1{F^2}
  u^3\big)u_r+
  r \frac 1{\bar C} u^3 \big[\frac 1r(r u_r)_r\big]_r\big\}_r=0.
   \end{matrix}
   \ee
   Here $E$ is again the evaporation number, $F$ is the Froude
   number, and $\e= \frac {h_0}L$ is a small parameter.
   Observe a rather complicated combination of various absorption and reaction-like
  non-divergent terms (with different nonlinear powers $u^3$, $u^4$, and $u^7$) in the first line of equation
  (\ref{oron2}).
Various exact solutions of non-conservative TFEs  can be found in
 \cite[Ch.~3]{GSVR}, where more references and a survey on TFE
theory are given.



Modern nonlinear parabolic theory and application to thin film
models demand better understanding of interaction of various
nonlinear terms and operators of different orders that can create
rather complicated spatio-temporal patterns and dissipative
structures. We chose one particular but special
 case of centre subspace behaviour that will be shown
to have rather robust mathematical significance.


\subsection{Basic limit model: the TFE with absorption}

 We study the large-time asymptotic behaviour of
nonnegative solutions of the {\em thin film equation} (TFE) {\em
with absorption} (for convenience, it is written for solutions of
changing sign to be studied also)
  \be
 \label{GPP}
 u_t = -\nabla \cdot (|u|^n \nabla \Delta u) - |u|^{p-1}u, 
   \ee
 where $n>0$ and $p>1$ are fixed exponents.
 Here we use the simplest second term which is not a
 differential operator but is represented by just a power
 function. Our main goal is to justify that in the critical case
 \be
\label{cr222} \textstyle{p_0= 1+n + \frac 4 N}
 \ee
various solutions of (\ref{GPP}) exhibit a complicated asymptotic
behaviour with some logarithmic corrections $\ln t$ for $t \gg1$.

We have chosen the  non-conservative equation (\ref{GPP}) for
simplicity and for better presentation of our mathematical tools.
 We claim that similar phenomena are quite general and appear also in
 various conservative models.
 Actually, the logarithmic correction $\sim (\ln \frac 1t)^{-1/7}$ in the behaviour for
 large enough $t $ was
rigorously observed \cite{Giac02} for the relaxed conservative
thin film model consisting of two operators,
 \be
 \label{rel1}
 u_t+(u^3 u_{xxx})_x +(u^n u_{xxx})_x=0, \quad \mbox{with}
 \quad 0<n<3 \quad (u \ge 0),
 \ee
 where the first term with
  $u^3$ corresponds to
 Reynolds'equation from lubrication theory.
 It was shown that, for concentrated enough initial data, in a certain intermediate time-range,
 the propagation rate is as follows:
  \be
  \label{Gia1}
   \mbox{$
  {\rm meas}\,\{u(x,t)>1\} \sim \big(\frac t{\ln \frac
  1t}\big)^{\frac 17},
   $}
   \ee
where the usual scaling-invariant factor $t^{\frac 17}$ is
associated with a standard dimensional analysis. Here,
 the log-correction is  a result of a delicate interaction of two
 scaling invariant operators in (\ref{rel1}). We believe that (\ref{Gia1}),
 proved in \cite{Giac02} rigorously, can be put into a framework
 of a centre manifold calculus (though a justification
 can be extremely hard).

 Log-corrections were observed for the limit stable
 Cahn--Hilliard equation \cite[Sect.~5.4]{EGW}
 \be
 \label{bb1}
  \mbox{$
 u_t = -\D^2 u + \D (|u|^{p-1} u), \quad \mbox{with}
 \quad p=1+ \frac 2N.
  $}
  \ee
For the semilinear case $n=0$ in the TFE (\ref{GPP}), such
logarithmically perturbed asymptotic are also well known and admit
 a rigorous mathematical treatment, \cite{GalCr}.


Thus,  we  consider for (\ref{GPP}) the standard free-boundary
problem (FBP) with {\em zero-height}, {\em zero contact angle},
and {\em zero-flux} ({conservation of mass}) conditions
 \be
 \label{GPP1}
 u=\n u=  {\bf \nu} \cdot (u^n \n \D u)  =0
 \ee
 at the singularity surface (interface) $\Gamma_0[u]$, which is
the lateral boundary of $ {\rm supp} \,u$ with the outward unit
normal ${\bf \nu}$.  Bounded, smooth and compactly supported
initial data
  \be
 \label{u00}
 u(x,0)=u_0(x) \quad \mbox{in} \,\,\,\Gamma_0[u]\cap\{t=0\}
  \ee
 are added to complete a suitable functional setting of the
 FBP.
 As usual, we assume that these three free-boundary conditions
  give a correctly specified problem for the fourth-order
parabolic equation, at least for sufficiently smooth and
bell-shaped initial data, e.g., in the radial setting.


Returning to basics of thin film theory, earlier references
  on derivation of the pure fourth-order TFE
  \be
 \label{TFE1}
 u_t= - \n \cdot (|u|^n \n \D u)
   \ee
and related models
   can be
found in \cite{Green78, Smyth88}, where first analysis of some
self-similar solutions  for $n=1$ was performed. Source-type  
similarity solutions of (\ref{TFE1}) for arbitrary $n$ were
studied in \cite{BPelW92} for $N=1$ and in \cite{BFer97} for the
equation in $\ren$. More information on similarity and other
solutions can be found in \cite{BerHQ00, BerHK00, Bow01}.
In general, the  TFEs are known to admit non-negative solutions
constructed by special ``singular" parabolic approximations of the
degenerate nonlinear coefficients; see the pioneering paper
\cite{BF1}, various extensions in \cite{Gr95, Ell96, EllS, LPugh,
WitBerBer} and the references therein. In what follows we study
the asymptotic behaviour  of sufficiently ``strong" weak solutions
of the TFEs, which satisfy necessary regularity and other
assumptions; see also the survey paper \cite{Beck05}. Notice that
 regularity theory for the TFEs is not fully developed,
especially in the non-radial $N$-dimensional geometry and for
solutions of changing sign, so we will need to impose extra formal
requirements, which are necessary for justifying our asymptotic
approaches.

Let us mention other  well-established and related conservative
thin film models with extra lower-order terms  describing the
dynamics of thin films of viscous fluids in the presence of two
competing forces; see \cite{BerPugh1}. For $N=1$, typical {\em
quasilinear} TFEs are 
 \be
\label{Bert1}
 u_t = -(u u_{xxx} + u^3 u_x)_x \quad (u \ge 0),
   \ee
and the general equation with power nonlinearities is
  \be
 \label{GEnm}
 u_t = - (u^n u_{xxx})_x - (u^m u_x)_x \quad (u \ge 0).
   \ee
We refer to papers \cite{Gl4, GBl6} and the book
\cite[Ch.~3]{GSVR} as sources of a large amount of further
references and results of TFE theory and application.


 In addition,
  our extra motivation of the TFEs model like
(\ref{GPP}) is mathematical and is associated with the previous
investigations of the quasilinear diffusion-absorption PDEs.

\subsection{A mathematical motivation: the PME with critical
absorption}

Second-order quasilinear parabolic equations with absorption are
well known in combustion theory. A key model is the {\em porous
medium equation} (PME) {\em with absorption}
  \be
  \label{PME1}
 u_t = \n \cdot (u^n \n u) - u^p  \quad \mbox{in} \,\,\, \ren \times
 \re_+ \quad (u \ge 0),
   \ee
where $n>0$ and $p$ are fixed exponents.
 A special interest to such equations was motivated by localised
 similarity solutions introduced by
L.K.~Martinson and K.B.~Pavlov at the beginning of the 1970s.
Mathematical theory of such PDEs was developed by A.S.~Kalashnikov
a few years later; see his survey \cite{Ka1} for the full history.
Besides new phenomena of localization and interface propagation,
for more than twenty years, the PME with absorption (\ref{PME1})
became a crucial model for determining various asymptotic
patterns, which can occur for large times $t \gg 1$ or close to
finite-time extinction as $t \to T^-$ (for $p<1$). For
(\ref{PME1}), there are a few
 parameter ranges with different asymptotics,
 \begin{gather*}
p>p_0=1+n + \textstyle{\frac 2N}, \quad p=p_0, \quad1+n < p <p_0, \quad p=1+n, \\
1< p < 1+n, \quad p=1, \quad 1-n < p<1, \quad p= 1-n, \quad p < 1-n,
\end{gather*}
 etc.; see references and details in \cite[Ch.~5,6]{AMGV}.

The most interesting and unusual transitional behaviour for
(\ref{PME1}) occurs at the first {\em critical} (or {\em Fujita})
absorption exponent
 \be
\label{cr11}
\textstyle{p_0= 1+n + \frac 2 N.}
 \ee
In this case (see details and references in \cite[p.~83]{AMGV}),
  the asymptotic behaviour as
$t \to \infty$ of nonnegative compactly supported solutions of
(\ref{PME1}) is described by the logarithmically perturbed
source-type solution of the pure PME,
  \be
 \label{PPME1}
\textstyle{u(x,t)  = (t \ln t)^{- \b N}[F(x/t^{\b}(\ln t)^{- \b
n/2}) + o(1)],\quad \mbox{where} \,\,\, \b = \frac 1{2+nN}.}
   \ee
  Without the logarithmic factors and the $o(1)$-term, the right-hand side is indeed
  the famous {\em Zel'dovich--Kompaneetz--Barenblatt} (ZKB) similarity
  source-type
  solution of the pure PME $u_t= \nabla \cdot (u^n \nabla u)$, which
  has the form
 \be
\label{sss1} \textstyle{u_s(x,t) = t^{- \b N} F(y), \,\,\, y =
x/t^\b , \quad \mbox{with} \,\,\, F(y)=\bigl[\frac {n \b}2
(a^2-|y|^2)_+ \bigr]^{\frac 1n},} \ee
  where $a>0$ is an arbitrary scaling parameter. This explicit
  solution
   dates back to the 1950. 
   In the class of
   solutions of changing sign, (\ref{PME1}) admits
a countable sequence of critical exponents, where the patterns
contain  similar logarithmic time-factors, \cite{GHCo}.

 \subsection{Outline of the paper: logarithmically perturbed patterns for the TFE with
 absorption}

\begin{figure}
\centering
\includegraphics[scale=0.85]{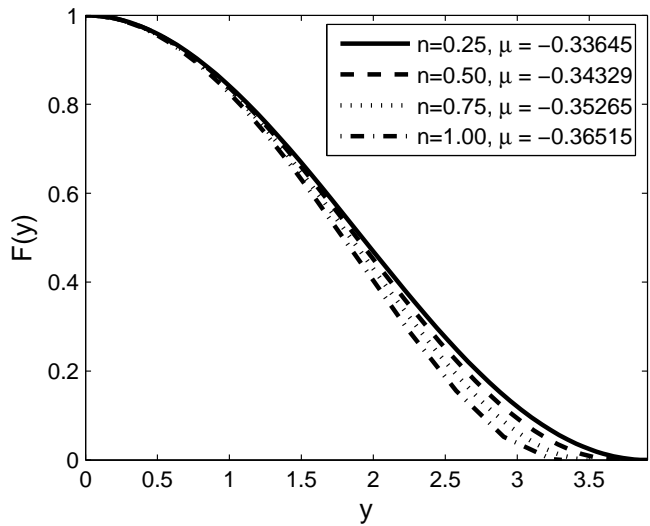}  
\vskip -.3cm
  \caption{The similarity profiles $F(y)$ of
(\ref{Od11}) for $N=1$ in four cases $n= \frac 14$, $\frac 12$,
$\frac 34$ and $1$; $F(0)=1$, $\mu= F''(0)$.} \label{Fig1}
\end{figure}

 In Sections \ref{Sect2} 
  we show that similar logarithmically perturbed
  source-type patterns exist for the TFE
 with absorption (\ref{GPP}), with the critical exponent
 (\ref{cr222}).
In this case, the source-type solutions of the TFE (\ref{TFE1})
take the form
 \be
\label{uss1} \textstyle{u_s(x,t) = t^{-  \b N} F(y), \quad y=
x/t^\b , \quad \mbox{with} \,\,\,\b = \frac 1{4+nN},}
  \ee
 where $F(y) \ge 0$ is a radially symmetric compactly supported solution
of the PDE \cite{BPelW92, BFer97}
 \be
\label{Od11} {\bf A}(F) \equiv -\n \cdot (F^n \n \D F) + \b \n F
\cdot y + \b N F=0.
  \ee
 In the case $n=1$, the similarity profile for the FBP
is given explicitly
 \be
\label{f112} \textstyle{F(y) = c_0 (a^2-|y|^2)^2, \quad c_0= \frac
1{8(N+2)(N+4)}, \quad a>0,}
  \ee
 and was first constructed in  \cite{Smyth88}.
Figure \ref{Fig1} shows  profiles $F(y)$ for $N=1$ in four cases
$n= \frac 14$, $\frac 12$, $\frac 34$ and $1$.
 The profiles are normalised by their values at $y=0$,
so $F(0)=1$.


First, for $n=1$, relying  on the explicit representation
(\ref{f112}) and good spectral properties of the corresponding
self-adjoint linearised rescaled operator,
  we show that, for  $p=p_0= 2+
\frac 4 N$, the TFE with absorption (\ref{GPP}) admits asymptotic
patterns of the following form:
 \be
\label{uln1} \textstyle{u(x,t) \sim (t \ln t)^{-\b N} F_*(x/t^\b
(\ln t)^{-\b N/4}) \quad \left(\b = \frac 1{4+N}\right).}
  \ee
  Here
$F_*$ is a fixed rescaled profile from the family (\ref{f112})
with a uniquely chosen parameter $a=a_*>0$ that depends on $N$
only. We also present evidence that similar logarithmic factors
can occur for arbitrary $n>0$ but this does not lead to
self-adjoint linearised operators and explicit mathematics. On the
other hand, for the semilinear case $n=0$, i.e., for the
fourth-order parabolic equation
  written for solutions
 of changing sign
 \be
\label{sem1} u_t = - \D^2 u -|u|^{p-1}u,
   \ee
   the critical
behaviour like (\ref{uln1}) is known to occur at the
 critical exponent $p= 1+ \frac 4N$ \cite{GalCr},
  which is precisely (\ref{cr222}) with $n=0$.
   In this case, the centre manifold analysis also uses spectral
    properties of a non self-adjoint linear operator studied
     in \cite[Sect.~2]{Eg4}.

In Section \ref{SectSu} we briefly describe the essence of the
easier supercritical case $p
> p_0$.  Very singular similarity solutions (VSSs)  in the
subcritical one $p \in (n+1,p_0)$ will be studied in a forthcoming
paper.


In Section \ref{Sect3}, we  explain how the critical asymptotic
behaviour occurs for the $2m$th-order TFE with absorption
 \be
\label{TFE2m} u_t = (-1)^{m+1} \n \cdot (u^n \n \D^{m-1}u) - u^p,
\quad m \ge 2,
 \ee
where the critical absorption exponent is
 \be
\label{crmmm} \textstyle{p_0= 1+n + \frac {2m} N,}
  \ee
  and again
$n=1$ leads to a simpler  self-adjoint case.

 In  Section
 \ref{Sect4}
   we discuss
similar local and global asymptotics for the Cauchy problem
admitting {\em maximal regularity} solutions of changing sign.

Finally, let us note that (\ref{cr222}) (and (\ref{crmmm}) for
equation (\ref{TFE2m})) is the critical {\em Fujita} exponent of
the TFE with source
 $$
 u_t = - \n \cdot (|u|^n \n \D u) + |u|^p \quad (n>0, \,\, p >1),
  $$
  i.e., for $p \in (n+1,p_0]$, all
  solutions with arbitrarily small initial data $u_0(x)$,
   where $\int u_0 >0$, blow-up in finite time, \cite{GPohTF}.

\section{\sc Rescaled equation and centre subspace behaviour}
\label{Sect2}


\subsection{To the style of the analysis}

 For
convenience of the Reader,
 we must emphasize   from the beginning that
all our final conclusions on centre subspace behaviours detected
below are {\em mathematically formal} when we deal with the
quasilinear case $n>0$. The semilinear  case $n=0$ is easier and
admits a rigorous treatment by invariant manifold theory,
\cite{GalCr}. It is then worth mentioning that
 there is no hope that such asymptotics can admit a reasonably
 simple rigorous treatment. We recall that even for the
 second-order model (\ref{PME1}) with $n>0$, there is no a full centre
 manifold justification of the main results that were proved by
 essential use of the Maximum Principle and comparison-barrier techniques;
 see \cite[Ch.~4]{AMGV}. Some of asymptotic patterns for (\ref{PME1})
  of centre subspace type  turned out to be very complicated,
  \cite{GHCo}. As we will show, the main difficulty is
 not a proper spectral theory of linearized operators
 (this is justified in many cases) but a
 justification of the centre subspace behaviour associated with
 such singular operators. On the other hand, we always clearly
 indicate the rigorous steps and split the whole approaches into a
 sequence of standard steps.
  We would be very pleased if some of our formal results and
  discussions would attract attention of experts in these areas of
  differential equations.

Thus, in what follows, we use by implication the following rule:

{\bf (i)} all conclusions concerning spectral and other properties
of self-adjoint singular elliptic and ordinary differential
operators are {\em rigorous} (or can be made rigorous after
sometimes technical manipulations; for non-self-adjoint cases we
are not that certain and extra analysis is necessary); and

{\bf (ii)} further extensions via above spectral properties to
describe the behaviour for TFEs close to center subspaces and
various matching procedures are mathematically {\em formal}.

\subsection{Rescaled equation}

We begin with rescaling the PDE (\ref{GPP}) with the critical
exponent (\ref{cr222}) according to the time-factors of the
source-type solution (\ref{uss1}), i.e.,  by setting
 \be
\label{2.1} u(x,t)= (1+t)^{-\b N} v(y,\t), \quad y= x/(1+t)^\b,
\quad \t = \ln (1+t),  \ee
  that leads to the following autonomous
rescaled equation in $\ren \times \re_+$:
 \be
\label{2.2}
v_\t = {\bf A}(v) - v^p,
 \ee
   where ${\bf A}$ is the operator specified in (\ref{Od11}). We
first need to check that a simple stabilization as $\t \to
+\infty$ to a nontrivial stationary solution in (\ref{2.2}) is not
possible.

 \begin{proposition}
\label{Pr.Non}
 The stationary equation
 \be
\label{St.1}
{\bf A}(g) - g^p = 0
 \ee
does not have a nontrivial compactly supported nonnegative
solution of the FBP.
\end{proposition}

 \noi {\em Proof.} Indeed, integrating (\ref{St.1}) over ${\rm supp} \, g$ yields
 $
\textstyle{\int g^p(y) \, {\mathrm d}y=0.} \quad \mbox{$\qed$}
$

\smallskip

This means that the only bounded nonnegative equilibrium for the
dynamical system (\ref{2.2}) is trivial,
 \be
 \label{equi1}
 g(y) \equiv 0 \quad \mbox{in} \,\,\, \ren.
  \ee

\smallskip

In order to detect the actual non-stationary asymptotic behaviour,
we next perform second rescaling by  introducing the as yet
unknown positive function $b(\t)$,
 \be
\label{2.3} v(y,\t)= b(\t)w(\zeta,\t), \quad \zeta = y/b^{\frac
n4}(\t),
 \ee
to get the following perturbed equation:
 \be
\label{2.4} \textstyle{ w_\t= {\bf A}(w) + \frac{b'}b \, {\bf C} w
- b^{p-1} w^p, \quad \mbox{where} \,\,\, {\bf C}w \equiv \frac n4
\, \n w \cdot \zeta - w.}
 \ee

\subsection{Linearisation}

Roughly speaking, in order to detect the asymptotic behaviour
according to (\ref{2.3}), we can use the estimate
 \be
\label{2.5}
b(\t) \approx \|v(\cdot,\t)\|_\infty \to 0 \quad \mbox{as} \,\,\, \t \to \infty,
 \ee
so that $\|w(\cdot,\t)\|_\infty \approx 1$ for $\t \gg1$. On the
other hand, in the radial setting, it is convenient to use $b(\t)$
for the scaling of the support of the solution $w(\zeta,\t)$ to
get that it approaches the unit ball $B_1$ as $\t \to \infty$; see
below.

We next perform the linearisation by setting
 \be
\label{2.6} w(\zeta,\t) = F(\zeta) + Y(\zeta,\t),
  \ee
  where $F$ is
a rescaled similarity profile from the family (\ref{f112}). Then
$Y$ solves the following  rescaled equation:
 \be
\label{2.7} \textstyle{Y_\t= {\bf A}'(F)Y + \frac{b'}b \, {\bf C}
F - b^{p-1} F^p +{\bf D}(Y) - b^{p-1}[(F+Y)^p-F^p],}
  \ee
  where
${\bf A}'(F)Y$ is the formal Frechet derivative of ${\bf A}$ at
$F$,
 \be
\label{2.8} {\bf A}'(F)Y=- \n \cdot(F^n \n \D Y) - \n \cdot (n
F^{n-1} Y \n \D F) + \b \n Y \cdot y + \b  N Y,
  \ee
  and ${\bf D}(Y)$
is a higher-order perturbation, which is quadratic in $Y \to 0$ on
smooth functions. Using the elliptic equation (\ref{Od11}) for
$F$,  on integration,
 \be
\label{2.9}
 F^n \n \D F= \b F y \quad \Longrightarrow \quad
  {\bf A}'(F)Y = - \n \cdot (F^n \n \D Y) + (1-n) \b  \n
\cdot (Y \zeta).
 \ee

\subsection{The self-adjoint case $n=1$}

It follows from (\ref{2.9}) that $n=1$ is a special case, where
the last term vanishes. We  fix $a=1$ in (\ref{f112}), so that
 the linearised operator  is 
 \be
\label{2.10} {\bf A}'(F)Y = - \n \cdot (F \n \D Y) \equiv - c_0 \n
\cdot ((1-|\zeta|^2) \n \D Y),
 \quad y \in B_1=\{|\zeta| < 1\}.
 \ee
One can see that it can be written in the form
 \be
\label{2.11} \textstyle{{\bf A}'(F)Y =- c_0 \frac 1
{\rho(|\zeta|)} \, [ \D(a(|\zeta|)\D Y) +2N \D Y], \quad
\mbox{where} \,\,\, a(|\zeta|) = 1-|\zeta|^2= \frac 1
{\rho(|\zeta|)},}
  \ee
  so, in the topology of $L^2_\rho(B_1)$,
operator (\ref{2.10}) is symmetric in $C^\infty_0(B_1)$
 with good coefficients, and hence
admits self-adjoint extensions. Next, using classical theory
\cite{BS}, we specify properties of its unique Friedrichs
self-adjoint extension. Its domain is constructed
by completing $C^\infty_0(B_1)$ in the norm induced by its positive 
quadratic form (corresponding to the operator $-{\bf A}'(F) - c_0
\D>0$)
\[
\textstyle{\langle Y,W \rangle_*  \equiv c_0 \int_{B_1} [a \D Y \D
W - (2N-1)  \n Y \cdot \n W].}
\]
The intersection of this Hilbert space with  the domain of the
maximal adjoint operator $D(({\bf A}'(F))^*) = \{v \in L^2_\rho:
\,\, {\bf A}'(F) v \in L^2_\rho\}$ defines  the domain of the
self-adjoint extension, which we denote by $D({\bf A}'(F))=
H^4_{\rho,0}$. In particular,  for any $v \in H^4_{\rho,0}$, there
holds
\[
\textstyle{v=0 \quad \mbox{on} \,\,\, \partial B_1,  \quad
\mbox{and} \quad \int_{B_1} a(\D v)^2 < \infty,}
\]
so that $H^4_{\rho,0} \subset H^2_{\rho,0}$.   Consider the
corresponding eigenvalue problem written in the form
 \be
\label{2.12} - c_0 [ \D(a(|\zeta|)\D \psi) + 2N \D \psi]= \rho \l
\psi \quad\mbox{in} \,\,\, H^4_{\rho,0}.
  \ee
 Since the embeddings
of the corresponding functional spaces $H^2_{a,0}$ and $H^1_0$
into  $L^2_\rho$ are compact, \cite[p.~63]{Maz}, we have that the
spectrum
 $\s({\bf A}'(F))$ is real and discrete.

For our purposes, it suffices to detect the eigenvalues and
eigenfunctions in the radial (ODE) setting  with the single
spatial variable $r=|\zeta| > 0$. The extension to the elliptic
setting is performed by using the polar coordinates $\zeta=(r,\s)$ in $B_1$,
  \be
\label{L-B} \D= \D_r + \mbox{$\frac 1{r^2}$} \,\D_\s,
  \ee
  where
$\D_\s$ is the Laplace-Beltrami operator on the unit sphere
$S^{N-1}= \partial B_1$ in $\ren$. $\D_\s$  is a regular operator
with a discrete spectrum in $L^2(S^{N-1})$ (each eigenvalue
repeated as many times as its multiplicity),
 \be
\label{LBsp} \s(-\D_\s) = \{\nu_k=k(k+N-2), \,\, k \ge 0\},
 \ee
and an orthonormal, complete,  closed subset $\{V_k(\s)\}$ of
eigenfunctions, which are homogeneous harmonic $k$-th order
polynomials restricted to $S^{N-1}$. We plug (\ref{L-B}) into
(\ref{2.11}),  where all the coefficients are radial functions,
and use the separation of variables
 \be
\label{2.15}
\psi(r,\s)= R(r) V_k(\s)
 \ee
for solving the eigenvalue problem (\ref{2.12}). For each
fixed $\nu_k$, we then arrive at a radial eigenvalue problem for $R$,
which is similar to that discussed below.

Thus we take $k=0$ in (\ref{2.15}) and consider the radially
symmetric eigenvalue problem (\ref{2.12}). For $N=1$, this problem
was studied in \cite{Bern02}, where further references are given.
It is not difficult to check that the radial operator ${\bf
A}'(F)$ has the discrete spectrum
 \be
\label{2.16} \s({\bf A}'(F))= \{ \l_k=c_0 k(k+2)(k+N)(k+N+2),
\,\,\, k=0,2,4,...\},  \ee
   where each eigenfunction $\psi$ is a
($k$+2)th-order polynomial,
 \be
\label{2.17} \psi_k(r) = b_k(r^{k+2}+...+d_k) \quad (\psi_k(1)=0),
 \ee where $\{b_k\}$ are normalization constants,  so that the
eigenfunction subset $\{\psi_k\}$ is orthonormal in $L^2_\rho$. In
particular,
 \be
\label{2.18} \textstyle{\psi_0(r)=b_0(r^2-1)>0, \quad b_0=- \sqrt{
\frac{N+2}{2 \o_N}} \quad (\l_0=0),}
  \ee
  where $\o_N=
\frac{2\pi^{N/2}}{N\Gamma(N/2)}$ is the volume of the unit ball in
$\ren$.  Such polynomials are complete and closed in typical
weighted $L^p$-spaces (a standard functional analysis result; see
\cite[Sect.~2.3]{Eg4} for details), and this justifies the
equality in (\ref{2.16}). Moreover, we then can use the
eigenfunction expansion with the orthonormal eigenfunctions subset
$\{\psi_k\}$ to deal with solutions of the corresponding PDE.

We next consider the rescaled equation (\ref{2.7}), which for
$n=1$ takes the form
 \be
\label{2.19} \textstyle{Y_\t= {\bf A}'(F)Y + \frac{b'}b \, {\bf C}
F - b^{p-1} F^p - \n \cdot (Y \n \D Y)- b^{p-1}[(F+Y)^p-F^p].}
  \ee
We deal with strong radially symmetric solutions of  (\ref{2.19}),
where we now choose the normalization function $b(\t)$ in
(\ref{2.3}) such that
 \be
 \label{supp1}
{\rm supp} \, w(\cdot, \t) = B_1 \quad \mbox{for} \,\,\, \t \gg 1.
 \ee
According to equation (\ref{2.19}), we then  need to assume that
$b(\t)$ is smooth, at least, for large $\t$, though this
requirement can be weaken by using a weak (integral) form of the
PDE. We now use the converging (in $L^2_\rho$ and in the
corresponding Sobolev class) eigenfunction expansion of the radial
solution
 \be
\label{2.20}
\textstyle{Y(\zeta,\t) = \sum_{k \ge 0} a_k(\t) \psi_k(\zeta)}
 \ee
to study the corresponding centre subspace behaviour for the
nonlinear operator ${\bf A}$. This part of our asymptotic analysis
is formal.

Thus substituting (\ref{2.20}) into (\ref{2.19}) and projecting
onto $\psi_0$ in $L^2_\rho$, we have that the first coefficient
satisfies the following perturbed ``ODE":
 \be
\label{2.21}
\textstyle{a_0'= - \g_1 \frac{b'}b - \g_2 b^{p-1} +... \,} ,
\,\,\, \mbox{where} \,\,\,
  \g_1 =- \langle {\bf C}F, \psi_0 \rangle _\rho > 0,
\,\, \g_2 = \langle F^p, \psi_0 \rangle _\rho >0.
  \ee
  We omit in
(\ref{2.21}) the higher-order terms assuming  that, for this type
of  behaviour, the non-autonomous perturbations are the leading
ones.  The signs of the coefficients $\g_{1,2}$ in (\ref{2.21})
are essential and are easily checked by integration.

It follows from (\ref{equi1}) and (\ref{supp1}) 
 that $b(\t) \to 0$ as $\t \to \infty$, so
 $$
 \mbox{$
 \frac {b'(\t)}{b(\t)} \quad \mbox{is
not integrable at $\t = \infty$}.
 $}
 $$
  Therefore, in order to have a
uniformly bounded expansion coefficient $a_0(\t)$, we need to
suppose that two terms on the right-hand side of (\ref{2.21})
annul each other asymptotically, so that, up to an integrable
perturbation, 
 \be
\label{2.23} \textstyle{\frac{b'}b=  - \frac{\g_2}{\g_1} \,
b^{p-1}+... \quad \mbox{for}\,\,\, \t \gg 1.}
  \ee
 This gives
the following necessary condition for existing of such a
behaviour:
 \be
\label{2.24} \textstyle{b(\t) = \g_* \t^{-\frac 1{p-1}} +... \, ,
\quad \mbox{where}\,\,\, \g_* =
\bigl[\frac{(p-1)\g_2}{\g_1}\bigr]^{-\frac 1{p-1}}.}
  \ee
   Returning to the original variables $\{x,t,u\}$, from
(\ref{2.24}) we obtain the asymptotic pattern (\ref{uln1}). The
rescaled profile $F_*$  is uniquely determined from (\ref{f112})
with $ a_*=\g_*^{n/4}$.

\subsection{Arbitrary $n \in (0, \frac 32)$.}
This non self-adjoint case is more difficult.
 Consider
the linearised operator (\ref{2.9}) for $n \not = 1$,
 where $F>0$ is the radial solution of the ODE (\ref{Od11}) in $B_1$;
  see \cite{BFer97} for existence, uniqueness, and  asymptotics.
   Then, for  $n < \frac 32$ \cite{BFer97},
 \be
\label{2.25}
 F(\zeta) \sim (1-|\zeta|)^2 \quad \mbox{as} \,\,\,
|\zeta| \to 1^-.  \ee
    Notice that there exists a one-parameter
family of the solutions given by
 \be
\label{2.26} \textstyle{F_a(\zeta) = a^{\frac 4n} F(\frac \zeta
a), \quad a>0.}
  \ee
    Firstly,  we claim that, for $n\neq1$, operator
(\ref{2.9}) is not symmetric in $L^2_\rho$ for any positive weight
$\rho$ in $B_1$; see Appendix A.  Secondly, we have that 
 \be
\label{2.27} \textstyle{\psi_0(\zeta) = \frac {\mathrm d}{{\mathrm
d} a}  F_a(\zeta)|_{a=1} \equiv  \frac 4 n \, F - \n F \cdot
\zeta}
  \ee
  is a positive eigenfunction of (\ref{2.9}) corresponding to
$\l_0=0$.
Observe that, with respect to the regularity, this eigenfunction
well corresponds to that for $n=1$; cf. (\ref{2.18}). Moreover, it
follows that, close to the singular point $|\zeta|=1$, the radial
part of (\ref{2.9}) is governed by the singular (at $\partial
B_1$) higher-order operator
 \be
\label{2.28} L_4 Y = -(s^{2n} Y''')', \quad s = 1-|\zeta|,
  \ee
which is symmetric in a weighted $H^{-1}$ topology (but we need a
result in $L^2$).
Solving the problem $L_4 Y = g$ with natural conditions at the
point $s=1$, which is assumed to be regular, we obtain, up to
compact perturbations, that
 \be
\label{2.29} \textstyle{L_2 Y \equiv - Y'' \sim \int^s s^{-2n}
\int^s g \equiv L_* g  \quad \Longrightarrow \quad Y \sim L_2^{-1}
L_* g ,}
  \ee
   where $L_2^{-1}$ is a compact operator in $L^2$.
It is easy to check that the integral operator $L_*$ is bounded in
$L^2$ for
 \be
\label{2.30} \textstyle{n < \frac 34,}
 \ee
  and then $L_2^{-1}
L_*$ is compact in $L^2$ as the product of a compact and a bounded
operator. Therefore ${\bf A}'(F)$ has discrete spectrum in the
parameter range (\ref{2.30}). This is not an optimal result since,
as we have seen, the discreteness of the spectrum remains valid
for $n=1$. We use this analysis as a simple  illustration of the
fact that the spectrum is usually discrete in the non-symmetric
case.

Thus $0 \in \s({\bf A}'(F))$ is an isolated eigenvalue. There is a
numerical evidence that the spectrum is discrete for all $n \in
(0, \frac 32)$; see \cite{Bern02}, where, moreover, first six
eigenvalues turned out to be {\em real} for $N=1$.
 Possibly this might mean that in a special topology of sequences as $l^2$ (not related
 to any of $L^2_\rho$) the linearised operator can be treated as
 symmetric and self-adjoint; cf. an example in \cite{Eg4}. 
 For $n=0$ in any dimension $N \ge 1$, the whole
spectrum is proved to be real. We refer to \cite[Sect.~2]{Eg4},
where this and other $2m$th-order operators were studied in
$L^2_\rho(\ren)$, i.e., for the Cauchy (not a free-boundary)
problem.

The rest of our study is formal. Once in the radial setting there
exists the centre subspace of ${\bf A}'(F)$, we are looking for a
(formal) centre subspace patterns for  (\ref{2.7})
 \be
\label{2.31} Y(\zeta,\t) = a_0(\t) \psi_0(\zeta) + ... \, .
  \ee
We assume the centre subspace dominance in the behaviour, so, as
usual, other terms in this expansion are assumed to be negligible
for $\t \gg 1$. Substituting (\ref{2.31}) into (\ref{2.7}), we
next find the projection onto the corresponding adjoint
eigenfunction $\psi_0^*$. In general, such an analysis becomes
rigorous if we establish existence of complete, closed and
bi-orthonormal eigenfunction subsets $\{\psi_k\}$ and
$\{\psi_k^*\}$. This  is an open problem except the case $n=1$
above and $n=0$ studied in \cite{Eg4}. We do not deal with the
adjoint operator ${\bf A}'^*(F)$ in this formal asymptotic
analysis. The projection onto $\psi_0^*$ yields the perturbed ODE
(\ref{2.21}), where the same coefficients $\g_{1,2}$ are
determined via the standard dual $L^2$ product, where  $\psi_0$ is
replaced by $\psi_0^*$. This formally leads to the same
asymptotics  (\ref{2.24}).

\smallskip

\noi{\bf The range $n \in [\frac 32,3)$.} The centre subspace
analysis applies also for larger $n$'s. The asymptotics of
similarity profiles change at $n= \frac 32$, where, instead of
(\ref{2.25}),
  \be
 \label{32.1}
 \mbox{$
 F(\zeta)  \sim(1-|\zeta|)^2 \bigl[ \frac 34 \, \b|\ln(1-|\zeta|)|\bigr]^{\frac 23}
 \quad \mbox{as} \,\,\, |\zeta| \to 1;
  $}
  \ee
 see \cite{BFer97}. On the other hand, for $n \in (\frac 32,3)$,
   \be
  \label{32.2}
   F(\zeta)  \sim(1-|\zeta|)^{\frac 3 n}
 \quad \mbox{as} \,\,\, |\zeta| \to 1.
  \ee
This regularity is sufficient for determining the corresponding
eigenfunction and the logarithmic behaviour.

For $n \ge 3$, the zero contact angle FBP
 does not provide us with a proper interesting evolution; see
 \cite{BFer97}.

\section{On the supercritical parameter range $p>p_0$}
\label{SectSu}

\subsection{Exponentially perturbed dynamical system for $p>p_0$}

Let us explain what we  expect for $p>p_0$ in  (\ref{GPP}). In
terms of the rescaled function
 \be
\label{v1} u(x,t)=(1+t)^{-\frac N{4+nN}} v(y,\t), \quad \t =
\ln(1+t),
 \ee
the equation takes the form
 \be
\label{v2}
\textstyle{v_\t = - \n \cdot (v^n \n \D v) + \frac 1{4+nN}\,
y \cdot \n v + \frac N{4+nN} \, v - {\mathrm e}^{- \g \t} v^p,}
  \ee
 where $\g = \frac{N(p-p_0)}{4+nN} > 0$ if $p>p_0$.
 Therefore the absorption  term $-u^p$ in (\ref{GPP})
generates an exponentially small perturbation in the rescaled
equation (\ref{v2}). Hence one can expect the convergence as $ \t
\to \infty$ to the rescaled similarity profile $F$ in (\ref{uss1})
of the limit mass, though the passage to the limit in (\ref{v2})
generates a number of technical difficulties. Here (\ref{v2}) is
formally an exponentially small perturbation of the autonomous
rescaled TFE
 \be
\label{v4} \textstyle{v_\t = {\bf A}(v)\equiv  - \n \cdot (v^n \n
\D v) + \frac 1{4+nN} \, y \cdot \n v + \frac N{4+nN} \, v}.
 \ee
As usual, we gain an extra advantage in the case $n=1$.

\subsection{The gradient case $n=1$}

It is known that, for $n=N=1$, the rescaled TFE (\ref{v4}) is a
gradient system, \cite{CarrT02}. Let us construct an
``approximate" Lyapunov function for strong solutions of the FBP
in $\ren$. Namely, we write down (\ref{v2}) in the form
 \be
\label{vv11} \textstyle{v_\t =  \n \cdot \bigl[ v \n \bigl(-\D v +
\frac1{2(4+N)} \,|y|^2 \bigr)\bigr] +  {\mathrm e}^{- \g \t} v^p}
 \ee
and multiply in $L^2(\ren)$ by $(-\D_v)^{-1} v_\t$, where, by
definition,
\[
(-\D_v)^{-1} w=g \quad \mbox{if} \,\,\, \D_v g \equiv \n \cdot (v
\n g) = -w,
\]
and $g=0$ at the free boundary of $v$. Then integrating by parts
yields the identity
 \be
\label{v5} \textstyle{\int v| \n (-\D_v)^{-1} v_\t|^2 = \frac
{{\mathrm d}}{{\mathrm d} \t} \, \bigl[ - \frac 12 \, \int |\n
v|^2 - \frac 1{2(4+N)} \, \int v |y|^2 \bigr] + J,}  \ee
  where $J$ corresponds
to the exponentially small term,
 \be
 \label{ee1}
\textstyle{ J= {\mathrm e}^{-\g \t} \int v^p (-\D_v)^{-1} v_\t.}
 \ee
Integrating (\ref{v5}) over $(0,T)$ yields
 $$
  \mbox{$
\int_0^T \int v| \n (-\D_v)^{-1} v_\t|^2 + \frac 12 \, \int |\n
v(T)|^2 +\frac 1{2(4+N)} \, \int v(T) |y|^2 \bigr]  \le C+
\int_0^T J,
 $}
 $$
so that, if the exponential term (\ref{ee1})
 $J \in L^1(\re_+)$, this yields
 extra uniform estimates,
\[
\sqrt v\,  \n (-\D_v)^{-1} v_\t \in L^2(\re \times \re_+) \,\,\,
\mbox{and}\,\,\,
 \n v, \, \sqrt v|y|
 \in L^\infty(\re_+;L^2).
\]
 Note that, obviously, (\ref{v5}) does not imply existence of
 a Lyapunov function (the non-autonomous PDE (\ref{vv11}) is
 not a gradient system). Anyway,
 since (\ref{v5})  gives a rather strong
estimate of $v_\t$ for $\t \gg 1$, this makes it possible to pass
to the limit $\t \to \infty$ and establish stabilization to an
equilibrium point (see the technique in \cite[p.~116-117]{AMGV}),
which is unique by the obvious mass-monotonicity with time of the
solution.

The symmetry of the Frechet derivative (\ref{2.10})
 at $F$ looks like a certain ``remnant" of the fact that
the original PDE is a gradient system.

\section{\sc Centre subspace patterns for the $2m$th-order TFE}
\label{Sect3}

We consider the $2m$th-order TFE with absorption (\ref{TFE2m})
   with the critical absorption (Fujita)
exponent (\ref{crmmm}).
      The proper
setting of a standard ``zero contact angle" FBP for the TFE
includes $m$+1 free boundary conditions at the free boundary
$\Gamma_0=
\partial \O(t) \times \re_+$ ($\O(t)$ is the support of $u(\cdot,t)$ at time
$t>0$),
 \be
\label{3.2.1} \textstyle{u= \n u=...= \frac
{\partial^{m-1}u}{\partial {\bf \nu}^{m-1}}= {\bf \nu} \cdot
\nabla (u^n \D^{m-1}u)=0 ,}
 \ee
  where $\nu$ is the unit outward normal to $\partial \O(t)$ that
is assumed to be sufficiently smooth.

\subsection{Similarity solutions}

The similarity solutions of the pure TFE
 \be
\label{3.2.2} u_t = (-1)^{m+1}  \n \cdot (u^n \n \D^{m-1}u)  \ee
take the standard form (\ref{uss1}) with
 \be
\label{3.3} \textstyle{\b= \frac 1{2m +nN}.}
  \ee
  One can see that
the critical exponent (\ref{crmmm}) is precisely the one, for
which the PDE (\ref{TFE2m}) possesses the same group of scaling
transformation. Then the rescaled profile $F$ satisfies the radial
restriction of the $2m$th-order elliptic equation
 \be
\label{3.4} {\bf A}(F)= (-1)^{m+1}  \n \cdot (F^n \n \D^{m-1}F) +
\b \n F \cdot  y + \b N F=0.
  \ee
  It seems that, for any $m \ge 3$, the
questions of existence and uniqueness of a  solution $F(y) > 0$ in
$B_1$ remain open. It is clear that, for large $m$, a standard
approach to existence based on a multi-parametric  shooting leads
to a complicated geometric
analysis (though some general conclusions in this  
geometry are likely).  We expect that the approach based  on the
$n$-branching (or a continuous homotopy connection with $n=0$) via
the classical theory \cite{VaiTr} makes it possible to explain
properties solutions, at least, for small $n>0$ by branching from
the linear case $n=0$ (but, surely, a standard approach to smooth
branching does not apply).
 For the Cauchy problem,
the spectral and other properties of the corresponding linear
operator (\ref{3.4}) for $n=0$ are given in \cite{Eg4}, and can be
used to clarify the behaviour for small $n>0$. For the FBP
(\ref{3.2.1}), an extra analysis of the linearised elliptic PDE is
necessary.

As usual, the case $n=1$ provides us with the explicit solution.
Writing the ODE (\ref{3.4})  in the radial divergent form (here
$y$ is actually $|y|$)
\[
(y^{N-1} F(\D^{m-1}F)')'=(-1)^{m}\b(y^N F)',
\]
on integration we obtain $\D^{m-1}F= (-1)^m \frac 12 \, \b y^2$.
Integrating this linear ODE $2m$-2 times yields the positive
solution in $B_1$
 \be
\label{3.5} \textstyle{F(y)= c_0 (1-|y|^2)^m, \quad\mbox{where}
\,\,\, c_0= \frac 12 \, \frac {N !!}{(2m)!!(2m+N)!!}.}  \ee

\subsection{Linearised operator}

We next follow the same scheme of the asymptotic analysis  as in
Section \ref{Sect2}. Similar to (\ref{2.8}), we introduce the
linearised operator
 \begin{align}
{\bf A}'(F)Y & = (-1)^{m+1} \n \cdot(F^n \n
\D^{m-1}Y)\label{3.6}\\ &+ \,(-1)^{m+1} \n \cdot (nF^{n-1} Y \n
\D^{m-1}F)+ \b \n Y \cdot y + \b NY.\notag
\end{align}
Using the ODE (\ref{3.4}), we have that
\[
(-1)^{m+1} \n \cdot (n F^{n-1} \n \D^{m-1}F) = - \b n N, \,\,\,
(-1)^{m+1}  nF^{n-1}  \n \D^{m-1}F = - \b ny,
\]
so (\ref{3.6}) can be written in the form
 \be
\label{3.7} {\bf A}'(F)Y= (-1)^{m+1} \n \cdot (F^n \n \D^{m-1}Y)+
\b N(1-n) y \cdot \n Y + \b N(1-n)Y,
  \ee
 and we again observe that $n=1$
is a special case.

\subsection{The self-adjoint case $n=1$}

Plugging the profile (\ref{3.5}) into (\ref{3.7}) yields the
following symmetric form of the operator:
 \begin{align}
{\bf A}'(F)Y&= c_0(-1)^{m+1} \n \cdot ((1-|y|^2)^m \n
\D^{m-1}Y)\label{3.8} \\ &\equiv  \,c_0(-1)^{m+1} [D^m((1-|y|^2)
D^m Y) + m(m-1)N \D^{m-1}Y],\notag
\end{align}
where $D^m$ denotes $\D^{m/2}$ for $m$ even and $\n \D^{(m-1)/2}$ for $m$ odd.
For instance, for $N=1$ and $m=3$, we have
\[
{\bf A}'(F)Y= c_0(1-y^2)^2[((1-y^2)Y''')''' + 6 Y^{(4)}].
\]
Having the symmetric operator (\ref{3.8}) in $C^\infty_0$, we next
determine its self-adjoint extensions, \cite{BS}. In particular
there exists the extension with discrete spectrum and polynomial
eigenfunctions in the radial setting (the non-radial case is
covered by using the spherical polynomials as in (\ref{2.15})).
The eigenvalues $\l_k$ for the polynomials $\psi_k(y)$ given in
(\ref{2.17}) are calculated by using (\ref{3.8}),
  \be
\label{3.9} \l_k = -c_0(k+2)k...[k+2-2(m-2)]
(k+N+2)(k+N)...[k+N-2(m-2)]
  \ee
for $ k=2(m-3),2(m-2),...\, .$ Using the eigenfunction expansion
in terms of complete and closed subset of polynomials $\{\psi_k\}$
partially justifies the asymptotic centre subspace analysis of the
corresponding rescaled equation (\ref{2.7}), which yields the same
ODE (\ref{2.21}) and hence the asymptotics (\ref{2.24}). Here in
the critical case (\ref{crmmm}) we still have $\frac 1{p-1}=\b N$
with $\b$ given by (\ref{3.3}). Finally, we arrive at the
asymptotic pattern (\ref{uln1}), where $4$ is replaced by $2m$.

\subsection{The general case $n \not = 1$}

We do not have such a self-adjoint operator, but anyway, once
$F>0$ in $B_1$ is determined, we obtain the radial eigenfunction
$\psi_0$ for $\l_0=0$ from the scaling symmetry group (\ref{2.26})
(the exponent $\frac 4n$ is replaced by $\frac{2m}n$) of equation
(\ref{3.4}). We can also guarantee that (\ref{3.7}) has compact
resolvent provided that $n>0$ is not large, so $\l_0=0$ is an
isolated eigenvalue.  The rest of the centre subspace behaviour
via the expansion (\ref{2.31}) remains unchanged and leads to
similar logarithmically perturbed asymptotic patterns. A rigorous
justification is a hard open problem.

\section{\sc Logarithmically perturbed patterns in the Cauchy problem}
\label{Sect4}

The asymptotic behaviour and similarity solutions for the TFE
(\ref{TFE1}) or (\ref{TFE2m}) posed in the whole space $\ren
\times \re_+$ are  less studied in thee literature.  For $n \in
(0,\frac 32)$, in the Cauchy problem (CP), the solutions
exhibiting the ``maximal regularity" at the interfaces are
oscillatory and of changing sign. See  \cite{Gl4, GBl6} and the
book \cite[Ch.~1]{GSVR}
 for correct meaning of the CP for thin film equations and further examples.
 For such solutions,  we need to assume that $u^n$ in (\ref{3.2.1}) is replaced by
$|u|^n$. Therefore  from now on in all the expressions and
equations we use the convention that
 \begin{equation}
\label{Conv1}
 \begin{gathered}
 u^n,f^n,v^n,w^n,... \,\,\, \mbox{are replaced by}
 \,\,\, |u|^n, |f|^n,|v|^n,|w|^n, ... \,\,\, \mbox{and}, \\
 u^p,f^p,v^p,w^p,... \,\,\, \mbox{are replaced by}
 \,\,\, |u|^{p-1}u, |f|^{p-1}f,|v|^{p-1}v,|w|^{p-1}w,... \, .
\end{gathered}
\end{equation}
We must admit that solutions of changing sign are less relevant
for many known physical applications of TFEs. Nevertheless, for
general PDE theory, it is key and of principal importance  to
include the Cauchy problem and  to show that the basic techniques
developed above apply to these much more complicated oscillatory
solutions.

The idea of sign changing solutions of TFEs is straightforward.
Indeed, the oscillatory properties of such solutions are a
manifestation of the fact that TFEs (\ref{3.2.2}) are
``homotopic", i.e., can be continuously deformed (e.g., as $n \to
0$) via non-singular uniformly parabolic PDEs with analytic
coefficients (see details in \cite[Sect.~14]{GBl6}) to the linear
{\em poly-harmonic equation}
 \be
\label{3.10} u_t = (-1)^{m+1} \D^m u \quad \mbox{in} \,\,\, \ren
\times \re_+.
  \ee
    By  classical parabolic theory (see e.g. Eidel'man
\cite{EidSys}),
 given initial data $u_0 \in L^1$, there exists the unique solution
  of the Cauchy problem for (\ref{3.10})
defined by the convolution
 \be
\label{3.11} u(x,t)= b(x,t) * u_0, \quad b(x,t) = t^{-\frac
N{2m}}F(y), \,\,\, y = x/t^{\frac 1{2m}},
  \ee
    where $b(x,t)$ is the fundamental
solution of the operator $D_t -(-1)^{m+1} \D^m$. For any $m \ge
2$, the rescaled kernel $F= F(|y|)$ is oscillatory as $y \to
\infty$, so this  property of changing sign is inherited by $L^1$
solutions of (\ref{3.10}). Assuming  a continuous (homotopic)
deformation of a class of solutions of (\ref{TFE1}) as $n \to
0^+$, this confirms that the TFE admits oscillatory solutions of
changing sign at least for not that large $n>0$.  Continuity and
homotopy concepts are effective for treating the Cauchy problem
for higher-order TFEs; see other  examples in \cite{GBl6}.

Then the  source-type solutions of the TFE take the same form
(\ref{uss1}), where the radial function $F$ of changing sign
solves the ODE (\ref{Od11}) with the convention (\ref{Conv1}). We
begin with the linear case $n=0$, which by continuity is going to
describe some properties of source-type solutions for sufficiently
small $n>0$.

\subsection{Properties of the rescaled fundamental solution
for $n=0$}

The linear ODE
 \be
\label{3.12} \textstyle{{\bf A}(F) \equiv - \D^2 F + \frac 14 \,
\n F \cdot y + \frac N 4 \, F=0  \quad \mbox{in} \,\,\, \ren}
  \ee
 is precisely the elliptic equation for the rescaled kernel $F$
of the fundamental solution in (\ref{3.11}). Therefore the
similarity profile $ F(y)$ exists and is unique under the
assumption
 \be
\label{3.12.1} \textstyle{\int F(y) \, {\mathrm d}y=1}
  \ee
  (in
view of existence-uniqueness of the fundamental solution).







Let us next describe an important relation between similarity
profiles for the FBP and  the Cauchy problem.
 Without loss of generality, we consider the case $N=1$, where on
 integration once (\ref{3.12}) takes the form
   \be
  \label{On.1}
  \mbox{$
  F'''= \frac 14 \, F y.
  $}
   \ee
It is easy to find all  decaying profiles corresponding to the CP
with the exponential WKBJ asymptotics as $y \to +\infty$,
  \be
\label{fss}
 \mbox{$
F(y) \sim y^{-\frac 13} {\mathrm e}^{a y^{ 4/3}}, \quad \mbox{with
$a$ satisfying} \,\,\, a^3=  \frac 14 \, (\frac 34)^3. $}
  \ee
There exist two complex conjugate roots for exponentially decaying
profiles
   \be
  \label{a01}
  \mbox{$
 a_{\pm}= - \frac 38 \, 4^{-\frac 13}(1 \pm {\rm i}\sqrt 3) \equiv -c_1 \pm {\rm i}c_2.
 $}
  \ee
This yields a {\em two-dimensional} bundle of oscillatory
solutions with the behaviour
  \be
 \label{osc.1}
 F(y) \sim y^{-\frac 13} {\mathrm e}^{-c_1 y^{4/3}}\bigl[A_1 \cos\bigl(c_2
 y^{\frac 43}\bigr) + A_2 \sin \bigl(c_2 y^{\frac 43}\bigr)\bigr] \quad \mbox{as} \,\,\, y \to
 \infty,
   \ee
where $A_1$ and $A_2$ are arbitrary constants. The algebraic
factor $y^{-1/3}$ is obtained by a standard asymptotic WKBJ
method. We observe here the periodic behaviour with a {\em single}
fundamental frequency (a result we will refer to in the TFE
analysis below).

 \begin{proposition}
 \label{Pr.inf}
  For $N=1$,
 the rescaled profile of the Cauchy problem $F=F_\infty$ given by
  $(\ref{3.12})$, $(\ref{3.12.1})$ is the limit of FBP similarity
  profiles on bounded intervals,
   \be
  \label{osc.2NN}
  F_\infty = \lim F_k,
    \ee
where each $F_k(y)$ is defined on interval $(-y_k,y_k)$,
   \be
  \label{On.21}
 F_k(\pm y_k)= F_k'(\pm y_k)=0, \quad \mbox{and}
  \ee
    \be
   \label{On.2}
   \mbox{$
   y_k = \big(\frac \pi{c_2} \, k\big)^{\frac 34}(1+o(1)) \quad \mbox{as} \,\,\,
   k \to \infty.
   $}
     \ee
    \end{proposition}

\noi{\em Proof.} The geometric aspect of such a property is
obvious in view of the oscillatory behaviour in (\ref{osc.1}).
  The convergence as $k \to \infty$ follows from
straightforward computations related to the whole exponential
bundle including (\ref{osc.1}) and the growing counterpart
 $$
 \mbox{$
 F(y)=y^{-\frac 13} {\mathrm e}^{a_0 y^{4/3}}+... \, , \quad
 \mbox{with} \,\,\, a_0= \frac 34 \,4^{-\frac 13}.
 $}
  $$
Then solving the FBP problem (\ref{On.21}) yields the asymptotic
equality
 $
 \cos(c_2y_k^{4/3}+ {\rm const.})=0,
 $
 whence the asymptotics (\ref{On.2}). $\qed$

 \smallskip

We also expect the following {\em Sturm property} be valid:
  \be
 \label{On.5}
 F_k(y) \quad \mbox{has precisely $k$ zeros on $(0,y_k)$}.
   \ee
Such a zero-number property is easily seen for $k \gg 1$, but is
not obvious for smaller $k$'s.




\subsection{Similarity profiles for $n>0$: existence and uniqueness}

 \begin{proposition}
\label{Pr.Ex1} For $N=1$ and $n \in (0,1)$, the ODE
$(\ref{Od11})$,
 $(\ref{Conv1})$ in $\re$ admits a unique solution $F \in C^3$ of
unit mass. The solution $F(y)$ is symmetric, compactly supported
and is oscillatory near finite interfaces at $y= \pm y_0$.
 \end{proposition}


\noi{\em Proof.} For $N=1$ the ODE (\ref{Od11}) has the form
  \be
 \label{TF11}
 |F|^n F'''= \b F y, \quad y \in \re.
   \ee
  Dividing by $|F|^n$ and  setting $|F|^{-n}F=g$ yields
    \be
   \label{g11}
   \mbox{$
   (|g|^\a g)'''= \b g y, \quad y \in \re, \quad \a= \frac n{1-n}.
    $}
     \ee
   Then  existence and uniqueness of a compactly supported solution
  $F \in C^3$ for any $n \in (0,1)$ follows from the results in
  Bernis--McLeod
  \cite{BMcL91}.
 $\qed$

 \smallskip

  For $n \in [1,\frac 32)$  solutions of (\ref{TF11}) are less regular (see below), so
  the techniques in  \cite{BMcL91} do not apply directly, but we expect that the
  existence-uniqueness result remains valid and can be extended further to
  some interval $n \in[\frac 32,n_{\rm h})$; see below.

In Figure \ref{Fig3} we have shown these similarity profiles for
some $n >0$ including the linear case $n=0$ leading to the ODE
(\ref{On.1}) for the fundamental rescaled profile. Here we observe
convergence of the fundamental profiles as $n \to 0^+$, which is
justified rigorously if all the zeros are ``transversal" and
isolated except the last one; see
below.


\begin{figure}
\centering
\includegraphics[scale=0.8]{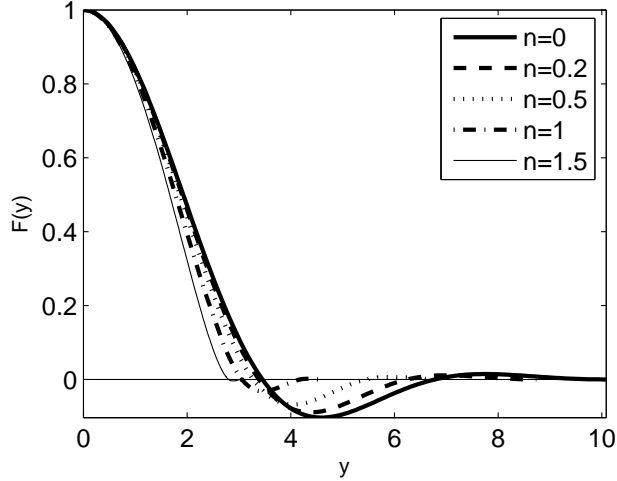}
\caption{The oscillatory CP profiles satisfying
\eqref{TF11}. Parameters of shooting are $F''(0)=-0.3379890$
($n=0$),
 $-0.3414702$ ($n=0.2$),
  $-0.3490986$ ($n=0.5$),
   $-0.3697143$ ($n=1$), and
    $-0.4052680$ ($n=1.5$).}
 \label{Fig3}
\end{figure}


\subsection{Oscillatory properties via periodic orbits}

We next describe the oscillatory properties of such changing sign
profiles $F(y)$ near interfaces. 
 We rescale $F$ to have
that
 $$
 {\rm supp} \, F= [-1,1].
 $$
  It was shown in
\cite{Gl4} that the asymptotic behaviour of $F(y)$ satisfying
(\ref{TF11}) near the interface point $y \to 1^-$ is given by the
expansion
 \be
\label{LC11}
 \mbox{$ F(y) = (1-y)^{\mu} \phi(s) , \quad s
=\ln(1-y), \quad \mu = \frac 3 n,
 $}
   \ee
where, after scaling $\phi \mapsto \b^{\frac 1 n} \phi$, the oscillatory component $\phi$ 
 satisfies the following autonomous ODE (we omit exponentially small terms):
 \be
\label{m=2.11} 
 \textstyle{\phi''' + 3(\mu-1) \phi''  + (3
\mu^2 - 6 \mu +2) \phi'+ \mu(\mu-1)(\mu-2)\phi +  \frac \phi
{|\phi|^{n}}=0.}
  \ee

\smallskip

\noi{\bf Oscillatory periodic orbits: existence.}   
  We
are now interested in  
 periodic solutions
$\phi_*(s)$ of  (\ref{m=2.11}), which according to (\ref{LC11}),
can determine the simplest typical  (and possibly stable and
generic) oscillatory behaviour of solutions near interfaces when
$s=\ln(1-y) \to - \infty$ as $y \to 1^-$. There are several
classic methods of ODE theory for establishing
 existence and multiplicity
of periodic solutions of finite-dimensional dynamical systems.
These are  various topological techniques, such as rotations of
vector fields, index, and degree theory; see
 \cite[Sect.~13,
14]{KrasZ}. 
 Another approach is based on
branching theory, 
\cite[Ch.~6]{VaiTr}. In our case, such an $n$-branching approach
is especially effective since for $n=0$ the unique
solution $F$ 
 is the rescaled kernel of the fundamental solution (a rigorous
justification of some aspects of branching for such degenerate
equations can be a hard problem). We also mention papers
 \cite{Ward, Liu, Kig} containing further related references
 and methods
  concerning modern theory of periodic solutions of
higher-order nonlinear ODEs. In general, 
equations like (\ref{m=2.11}) are a difficult object to study, and
especially the main difficulty is
 proving {\em uniqueness} of such periodic orbits.
 Therefore, later on, together with analytic techniques, we will need also to rely on
careful numerical evidence on existence, uniqueness, and stability
of periodic solutions.

It is curious that for $n=1$, the unique periodic solution can be
detected by a direct algebraic approach; see
\cite[Sect.~7.4]{Gl4}:

 \begin{proposition}
\label{Pr.n=1} For $n=1$, the ODE $(\ref{m=2.11})$ has a unique
$T$-periodic solution, with
  \be
 \label{per.1}
 T= -2 \ln s>\theta = 1.9248... \, , 
  \ee
 where $\theta=0.381966... $  is the unique root on the interval $(0,1)$ of the cubic equation
 \be
 \label{sys.3}
 \theta^3-2 \theta^2 - 2 \theta+1=0.
   \ee
 \end{proposition}

Indeed, for $n=1$, the nonlinearity in  (\ref{m=2.11}) is ${\rm
sign} \, \phi$ and the ODE
 is linear in the positivity and negativity domain
of solutions,
 $$
 \textstyle{\phi''' + 6 \phi''  + 11 \phi'+ 6\phi \pm 1 =0,}
 $$
so can be solved explicitly. Matching positive and negative
branches leads to the result.

Let us now state the main result concerning periodic orbits of the
ODE (\ref{m=2.11}).

\begin{theorem}
 \label{Th.Ex1}
 The ODE $(\ref{m=2.11})$ admits a nontrivial stable periodic solution
 $\phi_*(s)$ of changing sign for all
  \be
  \label{pp1}
   \mbox{$
  0<n<n_{\rm h} \in (\frac 32,n_{\rm +}), \quad \mbox{where}
   \quad n_{\rm +}= \frac 9{3+\sqrt 3}=1.9019238... \, .
    $}
    \ee
    \end{theorem}

    Uniqueness of such periodic $\phi_*(s)$ in the interval (\ref{pp1}) is still open.

\noi{\em Proof.} For the interval
 \be
 \label{in1}
 \mbox{$
0 < n < \frac 32,
 $}
 \ee
 the proof of existence is performed in \cite[p.~292]{Gl4}
 by a  shooting argument.
 Numerical representation of periodic solutions is given therein
 on p.~294; see also \cite[p.~143]{GSVR}.
  We need to point out the main two
 ingredients of the proof in \cite{Gl4}:

 (i) it is shown that for exponents (\ref{in1}) no orbits of the
 dynamical system (DS) (\ref{m=2.11}) are  attracted to infinity
as $s \to +\infty$, i.e., all orbits stay uniformly bounded; and

(ii) as a consequence, then (\ref{m=2.11}) is a dissipative DS
having a bounded absorbing set.

 Dissipative DSs are known to admit periodic solutions in rather
general setting \cite[Sect.~39]{KrasZ} provided these are
non-autonomous (so the  period is fixed). For the autonomous
system (\ref{m=2.11}), the proof in \cite[Sect.~7.1]{Gl4} was
completed by  shooting. Note that, in view of the last term,
(\ref{m=2.11}) is not a smooth dynamical system and solutions are
not locally $C^3$-smooth. Nevertheless, as shows local analysis
\cite[p.~291]{Gl4}, at least for $n \in (0,2)$, the nonlinearity
is integrable to guarantee local extensions of
 solutions through generic ``transversal" zeros. This means that
 the equivalent integral equation is well-posed and is composed
 from compact operators in a certain topology (this is necessary for application of classic methods
  of branching in Banach spaces, \cite[Ch.~7]{VaiTr}). We continue to
 deal with the differential equation, where the justification of
 calculus is done by local analysis.

It turns out that both properties (i) and (ii) also remain valid
for $n=\frac 32$,
so that a periodic solution $\phi_*$ also exists and is stable;
see Figure \ref{F32}. For the extension of $\phi_*$ to $n> \frac
32$, we will use the following crucial stability result:

\begin{figure}
 \centering
\includegraphics[scale=0.7]{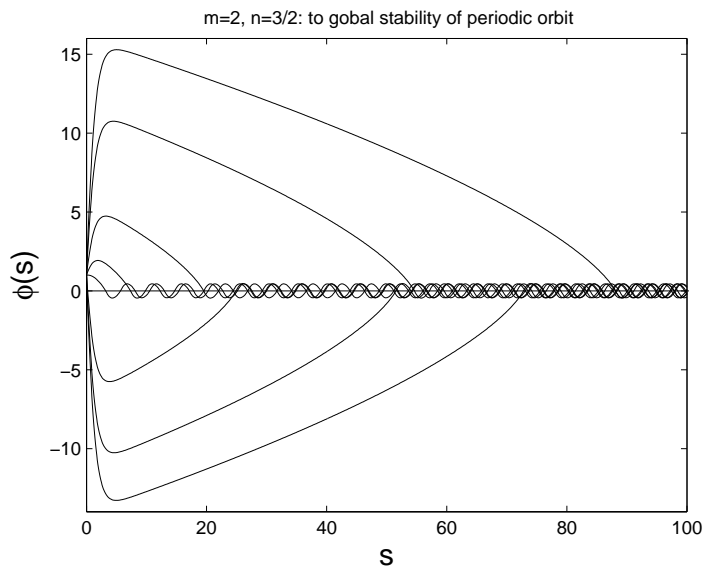}
\caption{\rm\small Convergence to the stable periodic solutions of
(\ref{m=2.11}) for $n= \frac 32$ for various Cauchy data posed at
$s=0$.}
 \label{F32}
\end{figure}

\begin{proposition}
\label{Pr.St1}
 If the periodic solution $\phi_*(s)$ of $(\ref{m=2.11})$
persists for all $\frac 32 \le n_{\rm h}<3$, then it is stable and
hyperbolic on this interval.
 \end{proposition}

 \noi{\em Proof.} Note that,
for $n \in (\frac32,3)$, there exist  two unstable
 constant equilibria of (\ref{m=2.11})  
  \be
 \label{eq.661}
 \mbox{$
 \phi_\pm =\pm \bigl[-\frac 1{\mu(\mu-1)(\mu-2)}\bigr]^{\frac 1n} \quad \mbox{for}
 \,\,\, n \in (\frac 32,3),
  $}
   \ee
   and we expect a stable periodic motion in between.
  Consider the eigenvalue problem for
  the ODE (\ref{m=2.11}) linearized about
 the $T$-periodic solution $\phi_*$ by setting $\phi=\phi_*+Y$,
  $$
  Y''' + 3(\mu-1)Y''+ (3 \mu^2-6\mu+2) Y'+
  \mu(\mu-1)(\mu-2)Y +(1-n)|\phi_*|^{-n} Y= \l Y.
  $$
 As usual, assuming that $\l \in {\mathbb C}$,  multiplying this by the complex conjugate
 ${\overline Y}$ in $L^2(0,T)$, taking the
 conjugate and multiplying by $Y$, and summing up both yields
  \be
  \label{ss1}
   \mbox{$
 -3(\mu-1) \int|Y'|^2+  \mu(\mu-1)(\mu-2)\int|Y|^2 + (1-n) \int|\phi_*|^{-n}
 |Y|^2= \frac {\l+\bar \l}2 \int |Y|^2.
 $}
  \ee
Since all the three terms on the left-hand side of (\ref{ss1}) are
negative for  any $\frac 32<n<3$, the result follows. The case
$n=\frac 32$ is similar since just the second term vanishes.
$\qed$

\smallskip

Thus, by classic branching theory, \cite[Ch.~6]{VaiTr}, stable
hyperbolic periodic solutions are locally extensible relative the
parameter $n \ge \frac 32$. In particular, using the hyperbolicity
of $\phi_*$ for $n= \frac 32$, we conclude that the periodic
solution exists in an interval $n \in [\frac 32, \frac 32+\d)$
with some $\d>0$, and the interval of existence must be open from
the right-hand side.

Finally, let us justify the estimate in (\ref{pp1}). To this end,
we multiply (\ref{m=2.11}) by $\phi'_*$ and integrate over $(0,T)$
to get for any $n \in (0,2)$
 $$
 \mbox{$
 - \int (\phi_*'')^2 + (3 \mu^2- 6 \mu +2) \int (\phi_*')^2=0,
  $}
  $$
  so that one needs
  $$
   \mbox{$
  3 \mu^2- 6 \mu +2>0 \quad \Longrightarrow \quad \mu = \frac 3n > \mu_+= \frac
  3{n_+}= \frac{3 + \sqrt 3}3.
   $}
   $$
 This completes the proof of Theorem \ref{Th.Ex1}. $\qed$

 \smallskip

\noi{\bf On heteroclinic bifurcation.} Since the periodic orbit
$\phi_*(s)$ remains stable and hyperbolic in the whole interval of
existence (\ref{pp1}), the end point $n=n_{\rm h}$ cannot be any
kind of subcritical saddle-node bifurcation, at which  two
branches meet each other. Classic bifurcation and branching theory
 \cite{KrasZ, VaiTr} then suggests that at $n=n_{\rm h}^-$ the DS
(\ref{m=2.11}) undergoes a {\em heteroclinic bifurcation} when the
period increases without bound (this claim needs further study and
a full analytical justification); see standard scenarios in
 Perko \cite[Ch.~4]{Perko}. Note that, by Proposition \ref{Pr.St1},
 the heteroclinic orbit occurred remains stable and hyperbolic.

\begin{figure}
\centering \subfigure[formation as $n \to  n_{\rm h}^-$]{
\includegraphics[scale=0.52]{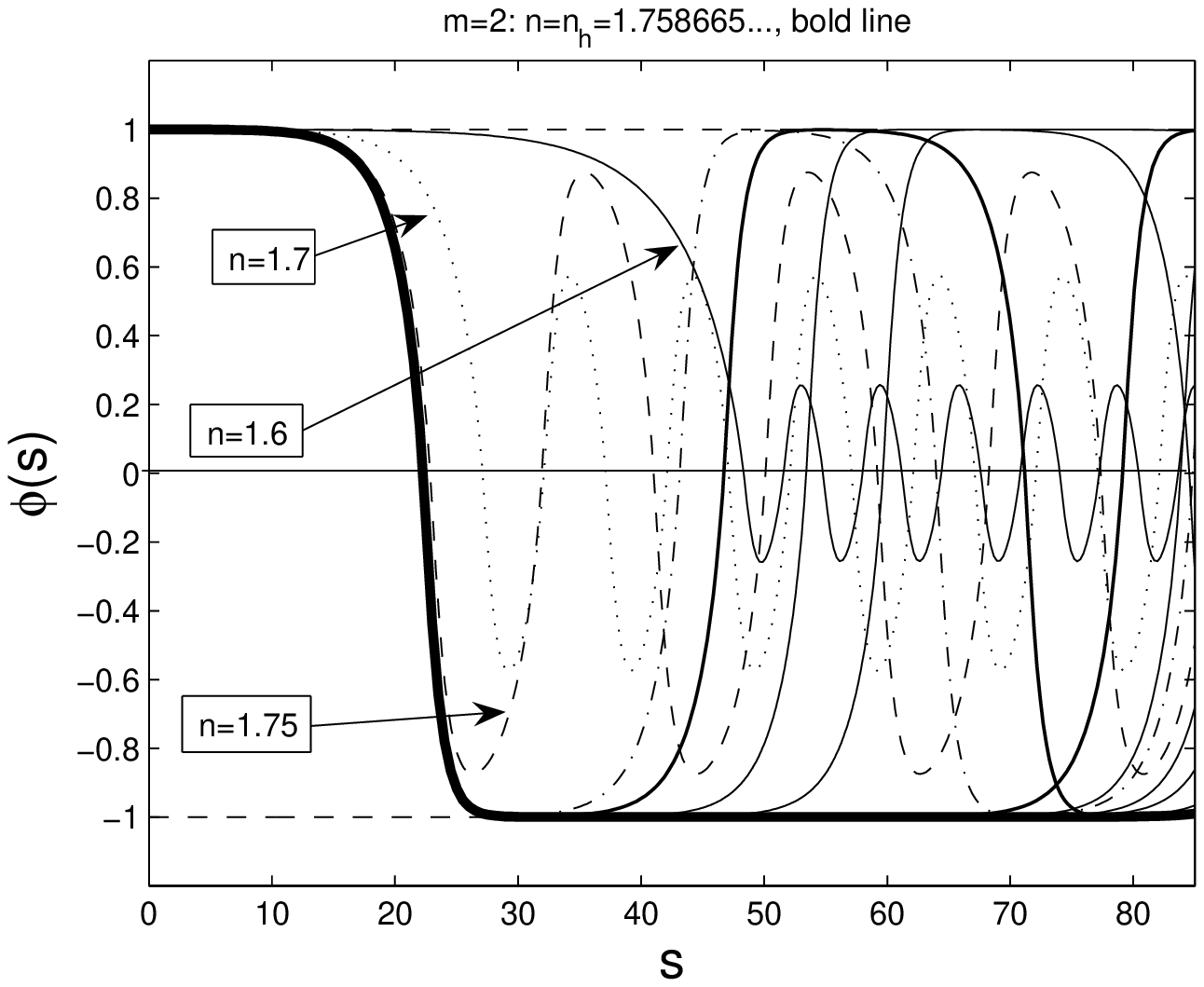}
} \subfigure[formation as $n \to  n_{\rm h}^+$]{
\includegraphics[scale=0.52]{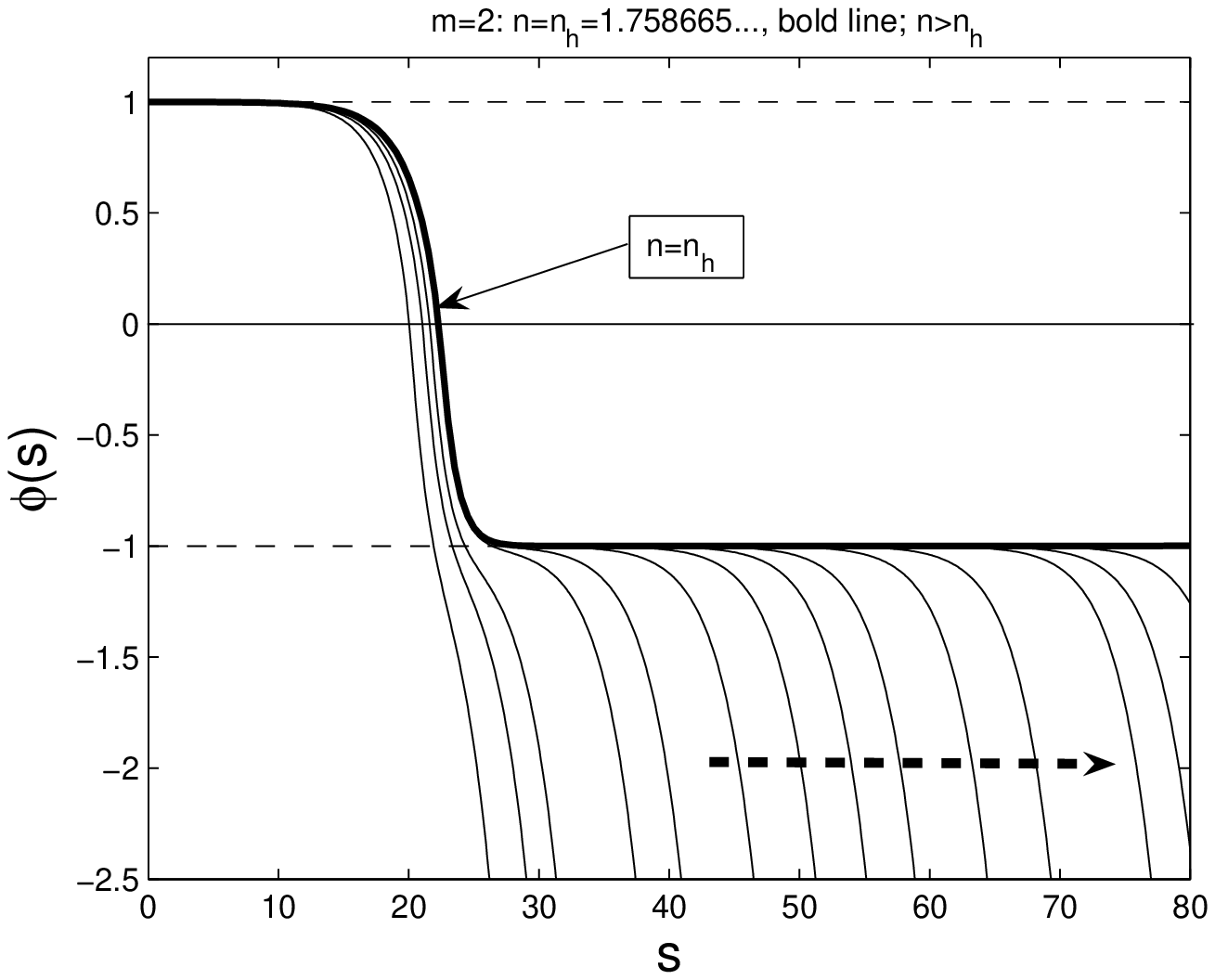}
}
 \vskip -.3cm
\caption{\rm\small Formation of a heteroclinic orbit as $n \to
n_{\rm h}$.}
 \label{FHet1}
\end{figure}

Numerically, $n_{\rm h}$ is given by
  \be
\label{n**1}  n_{\rm h}= 1.7598665026... \, .
  \ee
 Figure \ref{FHet1} shows formation of the heteroclinic orbit in
 both limits: as $n \to  n_{\rm h}^-$ (a) and $n \to  n_{\rm h}^+$
 (b).
This bifurcation
  exponent $n_{\rm h}$  
 plays the important
role and shows the parameter range of $n$'s, for which  many ODE
profiles near interfaces are  {\em oscillatory}
except those that approach the interface point $s=-\infty$ the stable manifold
of the constant equilibrium (\ref{eq.661}). In the interval (\ref{in1}), this manifold of
orbits of constant sign is empty,
so that all the orbits near $s=-\infty$ are oscillatory and coincide with
the periodic one $\phi_*(s+s_0)$, where $s_0 \in \re$ is a parameter of shifting.    
Indeed, this also characterizes important oscillatory features 
of the PDE. Note that some kind of a ``heteroclinic bifurcation"
phenomenon also exists for  the sixth-order ($m=3$) and
higher-order TFEs with more difficult mathematics involved; see
\cite[Sect.~13]{GBl6} and \cite[p.~142-147]{GSVR}.



\smallskip

\noi{\bf On 1D shooting for $n \in (1,n_{\rm h})$.} As a key
application of the above oscillation analysis, we have that
according to (\ref{LC11}), for all $n \in (0,n_{\rm h})$, there
exists a 1D bundle of oscillatory orbits of changing sign
 \be
 \label{ll1}
 F(y)=(1-y)^{\frac 3n}\phi_*(\ln(1-y)+s_0)+... \, ,
 \ee
 where $s_0 \in \re$ is an arbitrary parameter of phase shift in
 the periodic orbit $\phi_*(s)$. Recall that, for the ODE
 (\ref{TF11}), we need to shoot just a single symmetry condition
 at the origin,
  \be
  \label{ll2}
  F'(0)=0 \quad (F(0) \not = 0),
   \ee
so the 1D bundle (\ref{ll1}) is well-suited for this. In view of
oscillatory character of the behaviour in (\ref{ll1}), it is not a
great deal to prove the existence of such a $s_0$ to satisfy
(\ref{ll2}), while uniqueness (as expected) remains open.

\smallskip

\noi{\bf Further comments about $n_{\rm h}$.} For any $n>n_{\rm
h}$, the behaviour in the ODE (\ref{m=2.11}) becomes exponentially
unstable and we did not observe oscillatory or changing sign
patterns. This suggests that precisely above $n=n_{\rm h}$, the
ODE (and the corresponding PDE) loses its natural similarities
with the linear one for $n=0$ (though a continuous homotopic
connection is expected to be still available, {i.e.,} some local
properties of solutions dramatically change at $n_{\rm h}$).





Thus, in the range  $n \in (\frac 32,3)$, (\ref{m=2.11}) possesses
the positive constant solution $\phi(s) \equiv \phi_+$ given in
(\ref{eq.661}). This gives the behaviour (\ref{32.2}), so that,
for such solutions, formally, the FBP and the CP may coincide in
the ODE setting. But this is not the case for all the solutions
since for $n \in (\frac 32,n_{\rm h})$ there are other oscillatory
profiles with a similar (actually, a bit less) regularity at the
interfaces, so that the CP demands oscillatory solutions, while
the FBP can admit positive solutions; see more details in
\cite[Sect.~9]{Gl4}.
 In the parameter range $n \ge  n_{\rm h}$, 
the oscillatory behaviour is no longer  generic, so we  expect a
certain improvement of the positivity preserving properties of the
TFE, where the CP and the FBP may coincide; see further discussion
in \cite[Sect.~9.4]{Gl4}.
\subsection{The TFE with critical absorption $p=p_0$}

The formal asymptotics for the TFE (\ref{GPP}), (\ref{cr222}) is
 now calculated similarly using the centre subspace spanned by the eigenfunction (\ref{2.27}).
  Of course, we then do not gain any
explicit mathematics or symmetric operators as for $n=1$ in the
case of the FBP.

The main ideas of the analysis can be extended to the $2m$th-order
case, where many aspects of source-type and general solutions of
the Cauchy problem for the TFEs remain mathematically open. The
oscillatory character of solutions near the interface for $m=3$
was studied in \cite[Sect.~13]{GBl6}; see also
\cite[Sect.~3.7]{GSVR} for further examples for $m \ge 3$ and
other oscillatory PDEs.

\subsection{Supercritical range $p>p_0$}

We use the same scaling (\ref{v1}) and obtain the exponentially
perturbed rescaled PDE (\ref{v2}), which suggests that the
solutions behave as $t \to \infty$ as the source-type solution
with a finite positive mass attained at $\t=+\infty$ (no proof is
still available).

\smallskip
{\bf Acknowledgements.}  The authors would like to thank J.D.
Evans for discussions on thin film models with non-conservative
aspects, and A.~Leger for efficient consulting the authors with
numerical methods for higher-order ODEs.

 \begin{appendix}
\section{\sc The linearised operator is
not symmetric when $n \neq 1$ }
 \label{ApA}

  \begin{small}

We prove that, in the FBP setting, the linearised operator
 (\ref{2.9}) admits a  self-adjoint extension only when $n=1$.
Without loss of generality we consider the one-dimensional case,
and we formulate first the following results we are already
familiar with.

   \begin{proposition}
  The linearised operator $(\ref{2.9})$ in $\re$ is symmetric in some
weighted space $L^2_\rho$
   when $n=1$.
  \end{proposition}

  \noi {\em Proof.} For $N=n=1$, the linearised operator
   is given by
   \be
  \mbox{$
  {\bf A}'(f)Y=-(fY''')'-(Yf''')'+\frac 1 5 \,(Yy)'.
  $}
   \ee
  For this to be symmetric in $L^2_\rho$ with some weight $\rho \ge 0$, we require that
  \cite[Sect.~1]{Nai1}
   \be
  \label{SA}
  \mbox{$
  {\bf A}'(f)Y\equiv \frac 1 \rho \, \left[ (p_0 Y'')''-(p_1 Y')'+p_2 Y\right]
   $}
   \ee
  Expanding the right hand sides of these equations and comparing coefficients yields the following system:
   \begin{alignat}{2}
  &Y''''&&:\,\,\, \mbox{$ -f=\frac{p_0}{\rho}$}, \label{Y4}\\
  &Y'''&&:\,\,\, \mbox{$ -f'=\frac{2p_0'}{\rho}$},\label{Y3}\\
  &Y''&&:\,\,\, \mbox{$ 0=\frac{p_0''-p_1}{\rho}$}, \label{Y2}\\
  &Y'&&:\,\,\, \mbox{$ -f'''+\frac 1 5 y=-\frac{p_1'}{\rho}$},\label{Y1}\\
  &Y&&:\,\,\, \mbox{$ -f''''+\frac 1 5=\frac{p_2}{\rho}$}.\label{Y}
  \end{alignat}
  We know the exact solution of the ODE for $f$ when $n=1$ (see
  (\ref{f112})):
   \be
  \mbox{$
  f(y)=\frac{1}{120} \,(a^2-y^2)^2 \,\,\, \mbox{for}\: y\in(-a,a).
  $}
   \ee
  Substituting this into equation (\ref{Y}) yields $p_2=0$.
   Equation (\ref{Y1}) yields $p_1=C$ where $C$ is a constant. Equations (\ref{Y4})
    and (\ref{Y3}) yield $p_0^2=f$ and $\rho=-f^{-1/2}$.  Equation (\ref{Y2})
    is thus the consistency condition and is satisfied since it yields $p_1=C$
    (since $p_0''=p_1=C$).  Thus the
     linearised operator for the thin film equation is
      symmetric if $n=1$. $\qed$

   \begin{theorem}
  For $N=1$ and $n \neq 1$, operator $(\ref{2.9})$  is not symmetric in $L^2_\rho$
 for any weight $\rho>0$.
  \end{theorem}

  \noi {\em Proof.}
  The ODE for $f >0$ for any  $n > 0$ is
   \be
  \label{TFODE}
  \mbox{$
  -(f^nf''')'+\frac{1}{n+4} \,(fy)'=0.
  $}
   \ee
  The linearised operator (\ref{2.9}) is given by
   \be
  \mbox{$
  {\bf A}'(f)Y=-(f^nY''')'-n(f^{n-1}Yf''')'+\frac{1}{n+4} \, (Yy)'.
   $}
   \ee
  For this to be symmetric, we require identity (\ref{SA}) to hold.  Comparing coefficients yields
   \begin{alignat}{2}
  &Y''''&&:\,\,\, \mbox{$ -f^n=\frac{p_0}{\rho}$}, \label{Y4n}\\
  &Y'''&&:\,\,\,\mbox{$ -nf^{n-1}f'=\frac{2p_0'}{\rho}$}, \label{Y3n}\\
  &Y''&&:\,\,\, \mbox{$0=\frac{p_0''-p_1}{\rho}$},\label{Y2n}\\
  &Y'&&:\,\,\, \mbox{$-nf^{n-1}f'''+\frac{y}{n+4}=-\frac{p_1'}{\rho}$}, \label{Y1n}\\
  &Y&&:\,\,\,\mbox{$ -n(n-1)f^{n-2}f'''-nf^{n-1}f''''+\frac{1}{n+4}$}.\label{Yn}
  \end{alignat}
  From this
   \be
   \mbox{$
  p_0^2=f^n,\,\,\, p_1=p_0'',\,\,\, \rho=-f^{-n/2},\,\,\, p_2=\rho \bigl[-n(n-1)f^{n-2}f'''
  -nf^{n-1}f''''+ \frac 1{n+4} \bigr],
   $}
   \notag
   \ee
  and the consistency condition is
   \be
  \label{consistency}
  \mbox{$
  f^{\frac n2}(f^{\frac n2})'''=-nf^{n-1}f'''+\frac 1{n+4} \, y.
   $}
   \ee
To see if this coincides with equation (\ref{TFODE}) for some $f$
we use a Taylor expansion of $f(y)$ and check if
(\ref{consistency}) and (\ref{TFODE}) produce the same
coefficients for $f$.  To do this we set $f(0)=1,\,\,\,
f'(0)=f'''(0)=0$ and $f''(0)=b\in\mathbb{R}\setminus \{0\}$,
differentiate equations (\ref{consistency}) and (\ref{TFODE}) the
required number of times and set $y=0$.  The expansions coincide
up to the coefficient of $y^3$ but the coefficients of $y^4$ only
coincide if
   \be
  \label{b}
  \mbox{$
  b=\pm \frac{\sqrt{-6n(n^2+2n-8)(3n-2)}}{3n^3+6n^2-24n}.
  $}
   \ee
Since we require $b\in \mathbb{R}\setminus\{0\}$ we must have
$n\in (-4,0)\cup (\frac 23,2)$.  This gives us a range of values
of $n$, for which the linearised operator may be symmetric.  To
check whether it is we examine the coefficient of $y^6$ for
(\ref{consistency}) and (\ref{TFODE}).  If the operator is
symmetric, then the same value of $b$ should be obtained as in the
coefficients of $y^4$ for both equations.  For the coefficients of
$y^6$ to coincide we require
  \be
 \mbox{$
 b=0, \,\,\, \mbox{or} \,\,\,
  b=\pm
  \frac{2\sqrt{2}\sqrt{n(9n^3-40n^2-188n+464)(3n-2)}}{9n^4-40n^3-188n^2+464n}\,
  ,
   $}
  \ee
 and since we require
 $b\in \mathbb{R}\setminus\{0\}$ we discard $b=0$.
   For this $b$ to coincide with (\ref{b}) we require $n=\frac 23$.
    This contradicts the fact that we must have $n\in (-4,0)\cup (\frac 23,2)$ for
     the linearised operator (\ref{2.9}) to have a
      chance of being symmetric and admit a suitable (Friedrichs) self-adjoint
      extension.
        Hence the linearised operator is not symmetric if $n\neq 1$.
         $\qed$

\end{small}
\end{appendix}


\begin{thebibliography}{10}



   \bibitem
    {Beck05}
  J.~Becker and G.~Gr\"un, {\em The thin-film equation: recent advances and some new perspectives},
  {J.~Phys.: Condens. Matter}, {\bf 17} (2005), S291--S307.
 \bibitem 
  {Benney} D.J.~Benney, {\em Long waves on liquid films}, {J.~Math. and Phys.,} {\bf 45}
  (1966), 150--155.




\bibitem 
 {BF1}
 F.~Bernis and A.~Friedman, \emph{Higher order nonlinear degenerate
 parabolic equations}, J. Differ. Equat., \textbf{83} (1990), 179--206.


 \bibitem 
 {BerHK00}
 F.~Bernis, J.~Hulshof, and J.R.~King,
 {\em Dipoles and similarity solutions of the thin film equation in the half-line},
 { Nonlinearity,} {\bf 13} (2000), 413--439.


  \bibitem 
 {BerHQ00}
 F.~Bernis, J.~Hulshof, and F.~Quir\'os, {\em The ``linear" limit  of thin film flows
as an obstacle-type free boundary problem},
  {SIAM J.~Appl. Math.,} {\bf 61} (2000), 1062--1079.

\bibitem 
 {BMcL91}
 F.~Bernis and J.B.~McLeod, {\em Similarity solutions of a higher order nonlinear
 diffusion equation}, {Nonl. Anal.,} {\bf 17} (1991), 1039--1068.



\bibitem 
 {BPelW92}
 F.~Bernis, L.A.~Peletier, and S.M.~Williams, {\em Source type solutions
 of a fourth order nonlinear degenerate parabolic equation}, {Nonl. Anal.,}
   {\bf 18} (1992), 217--234.







  \bibitem
  {Bern02}
  A.J.~Bernoff and T.P.~Witelski, {\em Linear stability of
  source-type similarity solutions of the thin film equation},
  Appl. Math. Lett., {\bf 15} (2002), 599--606.


\bibitem{BerPugh1}
A.L. Bertozzi and M.C. Pugh, \emph{Long-wave instabilities and saturation in
  thin film equations}, Comm. Pure Appl. Math., \textbf{{LI}} (1998), 625--651.


  \bibitem 
{BS}  M.S.~Birman and M.Z.~Solomjak, {\rm Spectral Theory of
Self-Adjoint Operators in Hilbert Space}, D. Reidel,
Dordrecht/Tokyo, 1987.



\bibitem 
 {Bow01} M.~Bowen, J.~Hulshof, and J.R.~King, {\em Anomaluous exponents
 and dipole solutions for the thin film equation}, {SIAM J.~Appl.
 Math.,} {\bf 62} (2001), 149--179.


\bibitem
{CarrT02} J.A.~Carrillo and G.~Toscani, {\em Long-time asymptotic
behaviour for strong solutions of the thin film equations}, Comm.
Math. Phys., {\bf 225} (2002), 551--571.







  \bibitem{Eg4}
Yu.V.~Egorov, V.A.~Galaktionov, V.A.~Kondratiev, and
S.I.~Pohozaev,
  \emph{Asymptotic behaviour of global solutions to higher-order semilinear
  parabolic equations in the supercritical range}, Adv. Differ. Equat.,
{\bf 9} (2004), 1009--1038.




\bibitem{EidSys}
S.D. Eidelman, {Parabolic {S}ystems}, North-Holland Publ. Comp.,
  Amsterdam/London, 1969.



\bibitem 
{Ell96} C.M.~Elliott and H.~Garcke, {\em On the Cahn--Hilliard
equation with degenerate mobility}, {SIAM J.~Math. Anal.,} {\bf
27} (1996), 404--423.


\bibitem{EllS}
C.~Elliott and Z.~Songmu, \emph{On the {C}ahn-{H}illiard
equation}, Arch. Rat.
  Mech. Anal., \textbf{96} (1986), 339--357.


   \bibitem{Gl4}
J.D.~Evans, V.A.~Galaktionov, and J.R.~King, 
\emph{Source-type solutions of the fourth-order unstable thin film
equation}, Euro J.~Appl. Math., {\bf 18} (2007), 273--321.


      \bibitem
      {GBl6}   
J.D.~Evans, V.A.~Galaktionov, and J.R.~King, {\em Unstable
sixth-order thin film equation. I. Blow-up similarity solutions;
II. Global similarity patterns},
 {Nonlinearity}, {\bf 20} (2007), 1799--1841, 1843--1881.

 \bibitem
  {EGW}
J.D.~Evans, V.A.~Galaktionov, and J.F.~Williams, {\em Blow-up and
global asymptotics of  the limit unstable Cahn-Hilliard equation},
SIAM J.~Math. Anal., {\bf 38} (2006), 64--102.


\bibitem 
 {BFer97}
 R.~Ferreira and F.~Bernis, {\em Source-type solutions to thin-film equations in higher dimensions},
{European J.~Appl. Math.,} {\bf 8} (1997), 507--524.









  \bibitem
 {GalCr}
  V.A.~Galaktionov, {\em Critical global asymptotics in
  higher-order semilinear parabolic equations},
 Int. J. Math. Math. Sci., {\bf 60} (2003),  3809--3825.


\bibitem{GHCo}
V.A.~Galaktionov and P.J.~Harwin, {\em On evolution completeness
of nonlinear
 eigenfunctions for the porous medium equation in the whole space},
 Advances Differ. Equat., {\bf 10} (2005), 635--674.


\bibitem
 {GPohTF}
 V.A.~Galaktionov and S.I.~Pohozaev, {\em Blow-up and critical  exponents for  parabolic equations with
 non-divergent operators: dual porous medium and thin film
 operators}, J.~Evol. Equat, {\bf 6} (2006), 45--69.


\bibitem
 {GSVR} V.A.~Galaktionov and S.R.~Svirshchevskii, Exact Solution and
 Invariant Subspaces of Nonlinear Partial Differential Equations in Mechanics and Physics,
  Chapman$\,\&\,$Hall/CRC, Boca Raton,
Florida,
 2007. 




\bibitem
  {AMGV}
 V.A.~Galaktionov and J.L.~V\'azquez, {\rm A Stability Technique  for Evolution Partial Differential Equations.
 A Dynamical Systems Approach}, {\rm Progr. in Nonl. Differ.
 Equat. and their Appl.,} {\bf 56},
Birkh\"auser Boston, Inc., MA, 2004.





\bibitem
 {Giac02}
  L.~Giacomelli and F.~Otto, {\em Groplet spreading: intermediate
  scaling law by PDE methods}, {Comm. Pure Appl. Math.,} {\bf 55}
  (2002), 217--254.





\bibitem
 {Govor05}
L.V.~Govor, J.~Parisi, G.H.~Bauer, and G.~Reiter, {\em Instability
and droplet formation in evaporating thin films of a binary
solution}, {Phys. Rev. E}, {\bf 71},  051603
 (2005).


\bibitem 
{Green78} H.P.~Greenspan, {\em On the motion of a small viscous
droplet that wets a surface}, {J.~Fluid Mech.,} {\bf 84} (1978),
125--143.

\bibitem 
{Gr95} G.~Gr\"un, {\em Degenerate parabolic equations of fourth
order and a plasticity model with non-local hardening}, {Z.~Anal.
Anwendungen,} {\bf 14} (1995), 541--573.






\bibitem  
{Ka1}   A.S.~Kalashnikov, {\em Some problems of the qualitative
theory of second-order  nonlinear
 degenerate parabolic equations,}
{Russian Math. Surveys,} {\bf 42} (1987), 169--222.




\bibitem
 {Kig}
 I.T.~Kiguradze and T.~Kusano, {\em Periodic solutions of nonautonomous ordinary
differential equations
 of higher order}, Differ. Equat., {\bf 35} (1999), 71--77.




\bibitem{KrasZ}
M.A. Krasnosel'skii and P.P. Zabreiko, {Geometrical Methods of
Nonlinear
  Analysis}, Springer-Verlag, Berlin/Tokyo, 1984.



\bibitem
{LPugh} R.S. Laugesen and M.C. Pugh, {\em Energy levels of steady
states for thin-film-type equations}, J. Differ. Equat., {\bf 182}
(2002), 377--415.





    \bibitem
 {Liu}
 Z.~Liu and Y.~Mao, {\em Existence theorems for periodic solutions of higher
order
 nonlinear differential equations}, J.~Math. Anal. Appl., {\bf 216} (1997), 481--490.






\bibitem
 {Maz}
 V.G.~Maz'ja, {\rm Sobolev Spaces}, Springer-Verlag, Berlin/Tokyo,
 1985.


 \bibitem 
 {Nai1}
 M.A. Naimark, {\rm Linear Differential Operators}, 
  Frederick Ungar Publ. Co.,
 New York, 1968.


 \bibitem 
 {Oron00} A.~Oron, {\em Nonlinear dynamics of three-dimensional
 long-wave Marangoni instability
 in thin liquid films},
  {\em  Phys. Fluids,}
  {\bf 12} (2000), 1633--1645.

\bibitem 
 {Oron97} A.~Oron, S.H.~Davies, and S.G.~Bankoff, {\em Long-scale
 evolution of thin liquids films}, {Rev. Modern Phys.,}
  {\bf 69} (1997), 931--980.

\bibitem 
 {Oron02} A.~Oron and O.~Gottlied, {\em Nonlinear dynamics of
 temporally excited falling liquid films},
  {Phys. Fluids,}
  {\bf 14} (2002), 2622--2636.




\bibitem 
{Perko} L.~Perko, {\rm Differential Equations and Dynamical
Systems}, Springer-Verlag, New York, 1991.






\bibitem 
{Smyth88} N.F.~Smyth and J.M.~Hill, {\em High-order nonlinear
diffusion}, {IMA J.~Appl. Math.,} {\bf 40} (1988), 73--86.

\bibitem 
{VaiTr} M.A.~Vainberg and V.A.~Trenogin, {\rm Theory of Branching
of Solutions of Non-Linear Equations}, Noordhoff Int. Publ.,
Leiden, 1974.



\bibitem
 {Ward}
 J.R.~Ward, {\em Asymptotic conditions for periodic solutions of ordinary
 differential equations}, Proc. Amer. Math. Soc., {\bf 81} (1981), 415--420.



\bibitem{WitBerBer}
T.P.~Witelski, A.J.~Bernoff, and A.L.~Bertozzi, \emph{Blow-up and
dissipation
  in a critical-case unstable thin film equation},  Euro J.~Appl.
  Math., {\bf 15} (2004), 223--256.



\end{thebibliography}

\end{document}